\documentclass[11pt,reqno]{amsart}
\pdfoutput=1
\usepackage{tikz}
\usepackage{amsmath}
\usepackage{amssymb}
\usepackage{xcolor} 
\usepackage[letterpaper, left=2cm, right=2cm]{geometry} 
\usetikzlibrary{decorations.pathreplacing, fit} 
\usepackage{amsthm, amsmath,amsfonts,amssymb,euscript,hyperref,graphics,color,slashed}
\usepackage{graphicx}
\usepackage{mathrsfs}
\usepackage{comment}
\usepackage{import}
\usepackage{tikz}
\usepackage{latexsym}
\usepackage[makeroom]{cancel}
\usepackage[font=normalsize]{caption}
\usepackage{subfigure}
\usepackage[title]{appendix}
\usepackage{enumerate}
\usepackage{mathrsfs}
\usepackage{empheq}   
\usepackage{bm}
\usepackage{bookmark}
\usepackage{arydshln}
\usepackage{lineno,hyperref}
\usepackage{cite}
\usepackage{float}
\usepackage[utf8]{inputenc}
\usepackage{dashbox}
\usepackage{ulem}
\usepackage{lipsum}

\def\inte#1{
\displaystyle\mathop{#1\kern0pt}^\circ }

\def\virgp{\raise 2pt\hbox{,}}
\def\cdotpv{\raise 2pt\hbox{;}}

\def\C{\mathop{\mathbb C\kern 0pt}\nolimits}
\def\DD{\mathop{\mathbb D\kern 0pt}\nolimits}
\def\EE{\mathop{{\mathbb E \kern 0pt}}\nolimits}
\def\K{\mathop{\mathbb K\kern 0pt}\nolimits}
\def\N{\mathop{\mathbb N\kern 0pt}\nolimits}
\def\Q{\mathop{\mathbb Q\kern 0pt}\nolimits}
\def\R{\mathop{\mathbb R\kern 0pt}\nolimits}
\def\SS{\mathop{\mathbb S\kern 0pt}\nolimits}
\def\ZZ{\mathop{\mathbb Z\kern 0pt}\nolimits}
\def\TT{\mathop{\mathbb T\kern 0pt}\nolimits}
\def\P{\mathop{\mathbb P\kern 0pt}\nolimits}

\newcommand{\beq}{\begin{equation}}
\newcommand{\eeq}{\end{equation}}
\newcommand{\ben}{\begin{eqnarray}}
\newcommand{\een}{\end{eqnarray}}
\newcommand{\beno}{\begin{eqnarray*}}
\newcommand{\eeno}{\end{eqnarray*}}

\newtheorem{thm}{Theorem}[section]
\newtheorem{lem}{Lemma}[section]

\newtheorem*{Main Theorem}{Main Theorem}
\newtheorem{theorem}{Theorem}[section]

\newtheorem{remark}[theorem]{Remark}

\numberwithin{equation}{section}

\allowdisplaybreaks

\title[Morawetz energy estimate for wave equation]{A Morawetz type energy estimate for wave equation in $\R^2$ and application to elastic waves}
\author{Ningan Lai$^{\dag}$$^1$}

\author{Silu Yin$^{*}$$^2$}

\author{Yi Zhou$^{\#}$$^3$}

\date{}

\AtEndDocument{%
  \par
  \medskip
  \begin{tabular}{@{}l@{}}%
      \textsc{$^\dag$Department of Mathematics, Shaoxing University, China}\\
    \textit{E-mail address}: \texttt{$^1$ninganlai@usx.edu.cn}\\
    
    \textsc{$^*$Department of Mathematics, Hangzhou Normal University, China}\\
    \textit{E-mail address}: \texttt{$^2$yins@hznu.edu.cn}\\
    
    \textsc{$^\#$School of Mathematical Sciences, Fudan University, China}\\
    \textit{E-mail address}: \texttt{$^3$yizhou@fudan.edu.cn}\\
  \end{tabular}}

\begin{document}

\maketitle
\newcommand\blfootnote[1]{%
\begingroup
\renewcommand\thefootnote{}\footnote{#1}%
\addtocounter{footnote}{-1}%
\endgroup
}

\begin{abstract}

In this paper, we introduce a modified scaling Morawetz multiplier, which produces a weighted Morawetz type energy (non-negative) estimate for the inhomogeneous wave equation in $\R^2$. With this estimate in hand, an alternative proof of global existence for the Cauchy problem of quasilinear wave equation with small and compactly supported data is given. What is more, such weighted Morawetz type energy estimate also works for certain wave system with multiple speeds, which can be used to prove the global existence of some admissible harmonic elastic wave system in $\R^2$.

\end{abstract}


\tableofcontents


\setcounter{section}{0}
\numberwithin{equation}{section}

\section{Introduction}
We consider a quasilinear wave equation in $\mathbb{R}_+\times\mathbb{R}^2$
\begin{equation}\label{phieq}
  \partial_t^2 \phi-\Delta\phi=g^{ijk}\partial_i\phi \partial_j\partial_k \phi+h^{ijkl}\partial _i\phi\partial_j\phi\partial_k\partial_l\phi,
\end{equation}
where $g^{ijk}$, $h^{ijkl}$ are constant coefficients satisfying 
$$g^{ijk}=g^{ikj},\ h^{ijkl}=h^{ijlk}$$
with $i,j,k,l \in\{1,2\}$.
 Repeated indices obey the Einstein summation convention.
The initial data for $\phi$ are
 \begin{equation}\label{phidata}
 \phi(0,x)=\epsilon \phi_0(x),\quad \partial_t\phi(0,x)=\epsilon \phi_1(x),
 \end{equation}
 where $\phi_0,\phi_1$ are smooth real functions with compact support and  the positive parameter $\epsilon$ represents the smallness of the initial data.
 
Since the solutions to linear wave equations have $(1+t)^{-\frac{n-1}2}$ (in $n$-dimensional case) decay rate in time, the small data Cauchy problem of quasilinear wave equations is well-posed globally in time without any algebraic condition when $n\geq4$. While in three dimensions ($n=3$), the decay rate is critical to integrate with $t$. To obtain global existence with small initial data, faster time decay rate is needed such that it is integrable in time.  Christodoulou \cite{chris} and Klainerman \cite{klain} introduced a kind of nonlinear structure, which is the so-called null condition, to ensure the 
global existence when $n=3$. For $n=2$, as usual, 
we assume that $g^{ijk}$ and $h^{ijkl}$ in \eqref{phieq} satisfy the first and second null conditions, i.e., 
\begin{align}
g^{ijk}\omega_i\omega_j\omega_k=0,\label{fstnull}\\
h^{ijkl}\omega_i\omega_j\omega_k\omega_l=0,\label{secnull}
\end{align}
where $\omega_1=\frac{x_1}{r}$ and $\omega_2=\frac{x_2}{r}$ for $r=|x|$. Alinhac  solved this supercritical case and proved its global well-posedness in \cite{Alin01}.  For scalar quasilinear wave equations, he showed that
\begin{thm}\label{thm1}
If the null conditions \eqref{fstnull} and \eqref{secnull} are satisfied, then the Cauchy problem \eqref{phieq}-\eqref{phidata} with compactly supported small initial data admits a unique global smooth solution.
\end{thm}
The key point of Alinhac's method is the so-called ``ghost weight" when proceeding as in the energy estimate by 
using a weighted multiplier $e^p\partial_t\phi$. The ``ghost weight" $e^p$ is an appropriate ``bounded" weight such that the quadratic form of $\partial\phi$ obtained by integration by parts is non-negative.  As one of the detailed examples, this subtle weighted multiplier may be chosen as $e^{\arctan(r-t)}\partial_t\Gamma^\alpha\phi$. The condition of compact support in Theorem \ref{thm1} was removed by Cai-Lei-Masmoudi \cite{CL18} for wave equations with special nonlinearities. By a new class of $L^\infty-L^\infty$ estimate to the linear wave equations, Hou-Yin \cite{HY20} proved the global well-posedness of the general 2-D null form wave equations with non-compactly supported initial data.

In the work of Alinhac \cite{Alin01}, he used the full Lorentz invariance of the wave operator. In fact, Zha \cite{Zha} showed that the Lorentz boost operator is not necessary, based on the Klainerman-Sideris's inequality established in \cite{KlainSid,S1,S2}. Without Lorentz boost, Li \cite{Ld} improved the result to obtain the uniformly bounded higher-order energy, see also Dong-LeFloch-Lei \cite{DL21} in which the hyperboloidal foliation method was used. But it seems still hard to employ Alinhac's "ghost weight" method to wave systems with multiple speeds (not Lorentz invariant). In fact, the ``ghost weight" relies on the speed of propagation. A uniform ``weight" is not suitable for all waves with different propagation speeds.  

The goals of this paper are twofold,  first we offer an alternative proof of Theorem \ref{thm1}. For this end, we introduced an another weighted multiplier (improved scaling Morawetz type multiplier) to obtain a modified Morawetz type energy equality. The novelty is that we can show the non-negativity of the obtained Morawetz type energy, which can be also applied to wave systems with different propagation speeds, like elastic waves. Our second goal is to show the global existence of some admissible harmonic elastic waves in $\R^2$ by using the Morawetz type energy estimate mentioned above.

As introduced by John in \cite{John60, John66}, the motion of displacement for compressible harmonic, isotropic, hyperelastic materials satisfies a scalar wave equation. But it is interesting to see that {\it the admissible harmonic elastic waves} proposed by An-Chen-Yu \cite{ACY} obey a quasilinear wave system with different shear wave speed and pressure wave speed:
\begin{equation} 
    \partial_t^2 U-c_2^2\Delta U-\left(c_1^2-c_2^2\right)\nabla(\nabla\cdot U)=\text{grad}\ G(\partial U).
\end{equation}
Fortunately, the aforementioned Morawetz type energy method can be applied to the following admissible harmonic elastic waves satisfying the first and second null conditions:
\begin{equation}\label{1.6}
\begin{split}
    &\partial_t^2 U-c_2^2\Delta U-(c_1^2-c_2^2)\nabla(\nabla\cdot U)\\
   =&\nabla\left[g_1\left(\nabla\cdot U\right)\left(\nabla^\bot\cdot U\right)\right]+\nabla\left[h_1\left(\nabla\cdot U\right)^2\left(\nabla^\bot\cdot U\right)+h_2\left(\nabla\cdot U\right)\left(\nabla^\bot\cdot U\right)^2+h_3\left(\nabla^\bot\cdot U\right)^3\right].
    \end{split}
\end{equation}
We prove the following result. A detailed version is given in Theorem \ref{thm62} in Section \ref{sec6}.
\begin{thm}\label{thmhar}
The Cauchy problem of the admissible harmonic elastic wave system \eqref{1.6} with small smooth and compactly supported initial data admits a unique global classical solution. 
\end{thm}

The study of the long time behavior for the three-dimensional compressible elastic waves was pioneered by John \cite{Jo84}. He showed that the radially symmetric elastic waves will blow up in finite time. Almost global existence was obtained by John \cite{Jo88} and Klainerman-Sideris \cite{KlainSid}. With the null conditions, Agemi \cite{Ag} and Sideris \cite{S1,S2} proved the global existence of compressible elastic waves in 3D. For the incompressible case in 3D, the global existence was eatablished by Sideris-Thomases \cite{STh1,STh2}. Lei-Wang \cite{LW15} obtained the boundedness of the higher-order energy. If adding the viscosity, we refer to Chen-Zhang \cite{CZ06}, Cui-Hu \cite{CH22}, Hu-Meng-Zhang \cite{HM25}, Hu-Tu-Wen \cite{HT25}, Hu-Zhao \cite{HZ20}, Lei-Liu-Zhou \cite{LL08}, Lin-Zhang \cite{LZ08}, Qian-Zhang \cite{QZ}, Zhu \cite{Zhu} and the references therein.

The two-dimensional compressible elastic waves form a quasi-linear multiwave-speed system that is supercritical, the global existence of which is still open. One of the primary challenges lies in deriving a critical-type higher-order energy estimate:
\begin{equation}
    \frac{d}{dt}E_\kappa(t)\leq C\frac{(E_\kappa)^a}{1+t},
\end{equation}
where $a>1$. In the incompressible scenario, Lei \cite{lei} explored a strong null structure of the nonlinearity and established the small-data-global-existence result under the Lagrangian coordinates. Wang \cite{Wa17} gave an alternative proof of the global existence result under the Euler coordinates,
while Lei-Sideris-Zhou \cite{LS15} proved its almost global existence earlier. By vanishing viscosity limit, Cai-Lei-Lin-Masmoudi \cite{CL19} also showed the global well-posedness. For viscoelastic fluids, the global well-posedness was obtained by Lin-Liu-Zhang \cite{LL05} and Lei-Zhou \cite{LZ05}. For more related discussions,  we refer to Hu-Lin \cite{HL16}, Hu-Zhao \cite{HZ202d}, Lei-Liu-Zhou \cite{LL07} and the references therein. 

There are also some studies on general quasilinear wave systems with different propagation speeds in 2D. Hoshiga \cite{Ho06,Ho17} considered quasilinear wave systems with quadratic nonlinear terms satisfying the so-called Strong Null-conditions: Null-form and Nonresonance-form. We mention that the admissible harmonic elastic waves are not included in Hoshiga's cases.

\underline{\it Purpose of this work:}
We investigate the applications of Morawetz type multiplier in two-dimensional wave equations and systems. 

 {\it Morawetz-type multiplier}: For the $n$-dimensional linear wave equation $\Box\phi=0$ , in 1961, Morawetz used the "Friedrichs $abc$ method" to find a suitable multiplier in the form of $a\partial_t\phi+b\cdot\nabla\phi+c\phi$ to show the decay of solutions. 
   Specially, we call $\partial_r\phi+\frac{n-1}{2}\phi:=M\phi$ the Morawetz multiplier. This produces 
\begin{align}
& \frac{d}{dt}\int_{\mathbb{R}^n}\phi_tM\phi dx+\int_{\mathbb{R}^n}\Big(\frac{\Omega \phi}{r}\Big)^2dx+\underbrace{\frac{(n-1)(n-3)}{4}\int_{\mathbb{R}^n}\frac{\phi^2}{r^3}dx}_{S}=0,\ n\geq4,\\
&\frac{d}{dt}\int_{\mathbb{R}^3}\phi_tM\phi dx+\int_{\mathbb{R}^3}\Big(\frac{\Omega \phi}{r}\Big)^2dx+2\pi \phi^2(0,t)=0,\ n=3.
\end{align}
In dimension two ($n=2$), the integral term $S$ is singular. This indicates that the multiplier $M\phi$ is not available in the 2D case. According to the scaling invariance of the linear wave equation, by using scaling Morawetz  multiplier $t\partial_t\phi+r\partial_r\phi+\frac{n-1}{2}\phi:=MS\phi$, there holds the Morawetz estimate:
\begin{equation}\label{1.10}
    \int_0^\infty\int_{\mathbb{R}^n}\frac{\phi_t^2+|\nabla \phi|^2}{r}dxdt\lesssim \text{the initial energy},\ n\geq 3.
\end{equation}
The Morawetz estimate is also broken down for $n=2$. This is because the classical Hardy inequality is singular in 2D, which is critical to show \eqref{1.10}. 

In this paper, we use a modified scaling Morawetz  multiplier $$Z\Phi=(t+1)\partial_t\Phi+r\partial_r\Phi+\frac{n-1}{2}\Phi,$$ which produces the modified Morawetz dilation ``energy":
\begin{equation}\label{1.111}
    \frac12\int_{\mathbb{R}^n}\Big[(t+1)|\partial_t\Phi|^2+(t+1)|\nabla \Phi|^2+2r\partial_r\Phi\partial_t\Phi+(n-1)\Phi\partial_t\Phi\Big]dx.
\end{equation}
We show \eqref{1.111}  equals a weighted Morawetz type energy:
\begin{equation}
\begin{split}
  \tilde{M}_1(\Phi)=& \frac{3-n}{2}\int_{\mathbb{R}^n}\frac{t+r+1}{4}(\partial_t\Phi+\partial_r\Phi)^2dx+\frac{3-n}{2}\int_{\mathbb{R}^n}\frac{t-r+1}{4}(\partial_t\Phi-\partial_r\Phi)^2dx\\
   &+\int_{\mathbb{R}^n}\frac{t+1}{2r^2}(\Omega \Phi)^2dx+(n-1)\int_{\mathbb{R}^n}\frac{1}{8(t+1)}|\Phi|^2dx.
\end{split}
\end{equation}
And $\tilde{M}_1(\Phi)$ is not singular. Especially, the dimension $n=2$ is \underline{critical} to ensure $\tilde{M}_1(\Phi)$ is non-negative. With this observation, we offer an alternative proof of global existence for the Cauchy problem of quasilinear wave equations with small and compactly supported data in 2D. 

{\it Wave systems}: For the coupled systems of quasilinear wave equations 
\begin{equation}\label{1.13}
    \partial_t^2 u_i-c_i^2\Delta u_i=C_{ij}^{ab}(\partial u)\partial_a\partial_b u_j,\qquad i=1,\cdots,m,
\end{equation}
The global well-posedness for the two-dimensional Cauchy problem with small initial data is still unknown. The essential reason is that in this supercritical case, the $L^\infty$ norm of the solution to the linear equations has a lower decay rate in time‌. In detail,
the system \eqref{1.13} lacks the Lorentz invariance. Usually, we use the Klianerman-Sideris' inequality to compensate for this, 
\begin{equation}
    \begin{aligned}
    \underbrace{\sum_{|\alpha|=2}\sum_{|\beta|\leq\kappa-2}\|\langle c_it-r\rangle\partial^\alpha\Gamma^{\beta}u_i(t)\|_{L^2}}_{\chi_\kappa(u_i)}\lesssim& E_\kappa^{1/2}(u_i)+\sum_{|\alpha|\leq\kappa-2}\|(t+r)\Box_i\Gamma^\alpha u_i\|_{L^2}\\
    \lesssim& E_\kappa^{1/2}(u_i)+\underbrace{\sum_{|\alpha|\leq\kappa-2}\sum_{\beta+\gamma=\alpha}\|(t+r)C_{ij}^{ab}(\partial\Gamma^\beta u)\partial_a\partial_b\Gamma^\gamma u_j\|_{L^2}}_{K}.
\end{aligned}\end{equation}
But the $L^\infty$ norm of $\partial\Gamma^\beta u$ can only kill $t^{1/2}$-weight in term $K$. If we consider using the null conditions, it will also cause  complications. While the null condition can absorb the $t+r$ weight in $K$, it will lead to the loss of first-order derivatives with respect to time or space, resulting in like $\partial\Gamma^\alpha u$ changes to be $ \Gamma^{\alpha+1} u$. For single wave equation, we can use Hardy type inequality in Lemma \ref{lemhardy} to overcome this difficulty. But this way is failed for the coupled system \eqref{1.13} because of $$\Big\|\frac{\Gamma^{\alpha+1} u_k}{c_kt-r}\cdot\underbrace{(c_kt-r)\partial_a\partial_b u_j}_{\text{not } \chi_2(u_j)}\Big\|_{L^2}\ \text{for}\ j\neq k.$$

I‌n Section \ref{sec6}, with our Morawetz type energy, we build a new $L^2-L^2$ estimate in Lemma \ref{lem6.9} instead of the Klainerman-Sideris' inequality. With this in hand, we prove the global existence of some admissible harmonic elastic wave system in 2D.

Before giving our detailed proof for the above two theorems, we sketch the outline by the following flowchart first. 

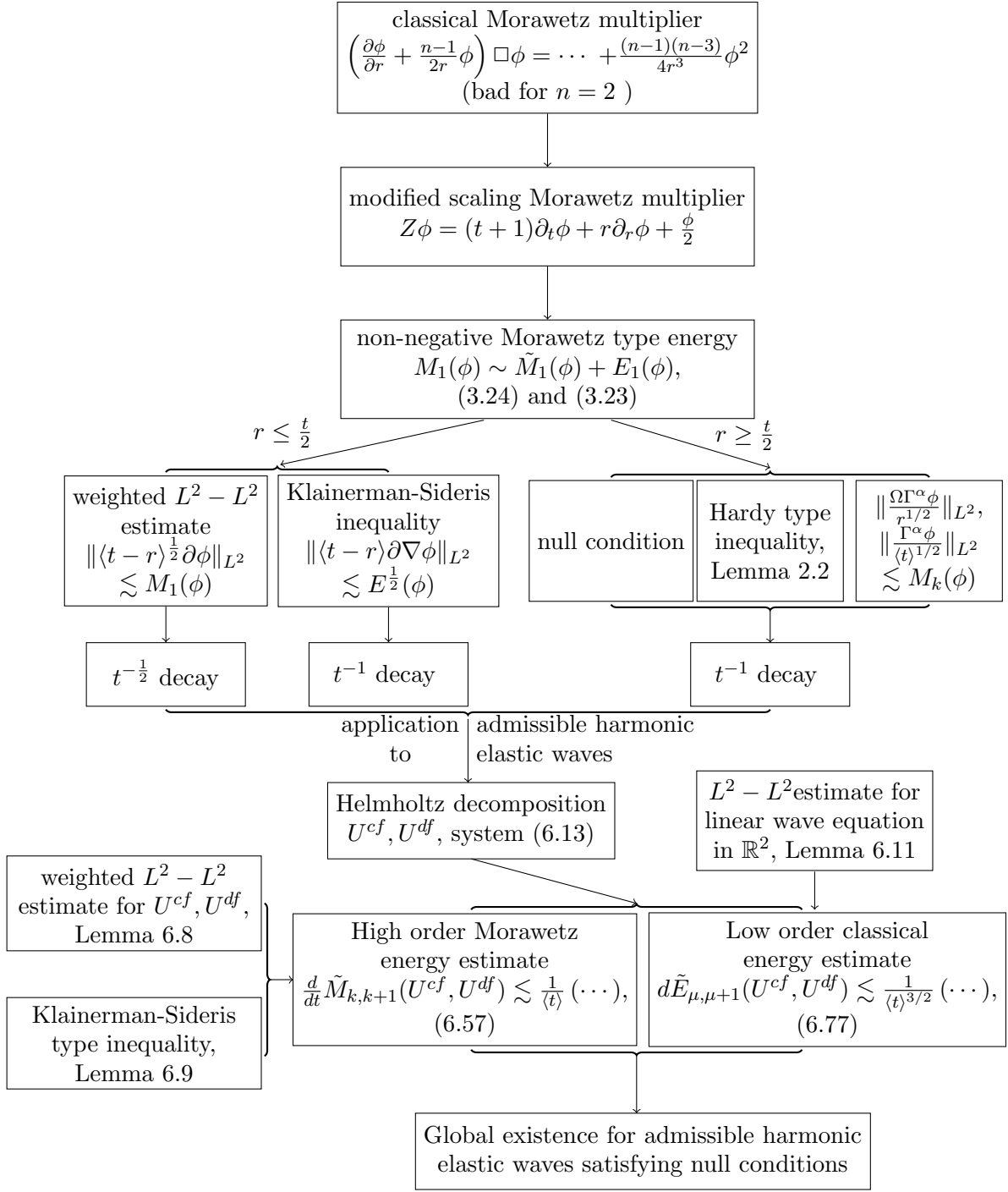
\begin{center}
\begin{tikzpicture}[    box/.style = {draw, rectangle, text width = 3.5cm, text height = 0.5cm, align = center},
    box1/.style = {draw, rectangle, minimum width = 6cm, minimum height = 1.55cm,  align = center},
    box2/.style = {draw, rectangle, minimum width = 2.3cm, minimum height = 1.8cm,  align = center},
    box3/.style = {draw, rectangle, minimum width =3.5cm, minimum height = 1.8cm,  align = center},
    box4/.style = {draw, rectangle, minimum width =2.5cm, minimum height = 1cm,  align = center},
    box5/.style = {draw, rectangle, minimum width =2.5cm, minimum height = 1.2cm,  align = center},
    box6/.style = {draw, rectangle, minimum width =2.5cm, minimum height = 1.2cm,  align = center},
    arrow/.style = {->, >=Stealth, thick},
    brace/.style = {decorate, decoration={brace, amplitude=5pt}, thick}
]

\node[box1] (11) at (-7,-0.5) {classical Morawetz multiplier \\ $ \left( \frac{\partial\phi }{\partial r} + \frac{n-1}{2r}\phi \right) \Box \phi= \cdots$ $ + \frac{(n-1)(n-3)}{4r^3} \phi^2 $ \\ \text{(bad for $n=2 $ )} };
\node[box1] (21) at (-7,-3) { modified scaling Morawetz multiplier \\$ Z\phi= \partial_t \phi + (t+1)\partial_t \phi  + r\partial_r  \phi + \frac{ \phi}{2}  $};
\node[box1] (31) at (-7,-5.4) { positive Morawetz type energy \\
$ M_1 (\phi) \sim \widetilde{M}_1 (\phi) + E_1(\phi)$,  \\ \eqref{eqvMM}  and \eqref{M1bar}};

\draw[->] (11.south) -- (21.north) node[midway, above] {};
\draw[->] (21.south) -- (31.north) node[midway, above] {};

\node[box2] (41) at (-13,-8.1) { weighted $ L^2-L^2 $ \\estimate \\$ \| \langle t-r\rangle^{\frac{1}{2}} \partial \phi \|_{L^2}$\\$ \lesssim M_1 (\phi) $ };
\node[box2] (42) at (-9.5,-8.1) { Klainerman-Sideris \\inequality \\$ \| \langle t-r\rangle \partial \nabla \phi \|_{L^2}$\\$ \lesssim E^{\frac{1}{2}} (\phi) $ };
\node[box2] (43) at (-6,-8.1) {null condition};
\node[box2] (44) at (-3.5,-8.1) { Hardy type \\inequality,\\ Lemma \ref{lemhardy} };
\node[box2] (45) at (-1,-8.1) {$ \| \frac{ \Omega \Gamma^\alpha \phi} {r^{1/2}} \|_{L^2}, $ \\
$ \| \frac{\Gamma^\alpha \phi} {\langle t\rangle^{1/2}} \|_{L^2} $
\\ $ \lesssim M_k (\phi)  $};

\draw[decorate, decoration={brace}, line width=0.8pt] (-13,-6.9) -- (-9.5,-6.9) node[midway, below=30pt] {};

    \coordinate (start) at (-8,-6.2);
    \coordinate (bend1) at (-11.2,-7);

    \draw[->, thin] (start) -- (bend1)node[pos=0.8, above left] {\textcolor{black}{$ r\leq \frac{t}{2}$}}; 

\draw[decorate, decoration={brace}, line width=0.8pt] (-6,-7.1) -- (-1,-7.1) node[midway, below=30pt] {};

    \coordinate (start2) at (-6,-6.2);
    \coordinate (bend2) at (-3.5,-6.95);

    \draw[->, thin] (start2) -- (bend2)node[pos=0.8, above left, xshift=20pt] {\textcolor{black}{$ r\geq \frac{t}{2}$}}; 

\node[box4] (51) at (-13,-10.2) { $ t^{-{\frac{1}{2}}} $ decay};
\node[box4] (52) at (-9.5,-10.2) { $ t^{-1} $ decay };
\node[box4] (53) at (-3.5,-10.2) { $ t^{-1} $ decay };

\draw[decorate, decoration={brace,mirror }, line width=0.8pt] (-6,-9.1) -- (-1,-9.1) node[midway, below=30pt] {};

\draw[->] (41.south) -- (51.north) node[midway, above] {};
\draw[->] (42.south) -- (52.north) node[midway, above] {};

    \coordinate (start3) at (-3.5,-9.2);
    \coordinate (bend3) at (-3.5,-9.7);

    \draw[->, thin] (start3) -- (bend3)node[pos=0.8, above left, xshift=20pt] {}; 

\draw[decorate, decoration={brace,mirror }, line width=0.8pt] (-13,-10.8) -- (-3.5,-10.8) node[midway, below=30pt] {};
    \coordinate (start4) at (-8.25,-10.9);
    \coordinate (bend4) at (-8.25,-11.9);

    \draw[->, thin] (start4) -- (bend4)
    node[pos=0.8, above left, xshift=0pt,align=center] {application \\to }  
    node[pos=0.8, above right, xshift=0pt, align=left]  {\ admissible  harmonic \\elastic waves}; 

\node[box6] (61) at (-8.2,-12.5) { Helmholtz decomposition \\ $ U^{cf}, U^{df}$,  system (6.13) };
\node[box6] (62) at (-2.8,-12.5) { $ L^2 - L^2  $estimate for\\ linear wave equation\\ in $\mathbb R^2$,  Lemma \ref{L2L}};

\node[box3] (71) at (-13.5,-13.8) {  weighted $ L^2-L^2 $ \\estimate  for $ U^{cf}, U^{df}$, \\ Lemma \ref{L68}};
\node[box2] (72) at (-8.2,-15) { High order Morawetz \\energy estimate\\ 
$ \frac{d}{dt} \widetilde{M}_{k,k+1} (U^{cf}, U^{df})\lesssim \frac{1}{\langle t\rangle} \left(\cdots \right) $, \\ \eqref{656} };
\node[box2] (73) at (-3,-15) {  Low order classical\\ energy estimate\\
$ d \widetilde{E}_{\mu,\mu+1} (U^{cf}, U^{df})\lesssim \frac{1}{\langle t\rangle^{3/2}} \left( \cdots \right) $, \\ \eqref{6.62}};
\node[box3] (74) at (-13.5,-16) { Klainerman-Sideris \\ type inequality, \\ Lemma \ref{lem6.9}};

\draw[decorate, decoration={brace}, line width=0.8pt] (-11.65,-13.8) -- (-11.65,-16.2) node[midway, below=30pt] {};

    \coordinate (start5) at (-11.6,-15);
    \coordinate (bend5) at (-10.6,-15);

    \draw[->, thin] (start5) -- (bend5)node[pos=0.8, above left, xshift=20pt] {}; 
    
\draw[decorate, decoration={brace}, line width=0.8pt] (-8.2,-14) -- (-3,-14) node[midway, below=30pt] {};

    \coordinate (start5) at (-8.2,-13.2);
    \coordinate (bend5) at (-5.6,-13.85);

    \draw[->, thin] (start5) -- (bend5)node[pos=0.8, above left, xshift=20pt] {}; 
    \coordinate (start6) at (-2.8,-13.2);
    \coordinate (bend6) at (-2.8,-14);

    \draw[->, thin] (start6) -- (bend6)node[pos=0.8, above left, xshift=20pt] {}; 

\node[box5] (81) at (-5.5,-17.7) { Global existence for admissible harmonic\\ elastic 
 waves satisfying null conditions};

\draw[decorate, decoration={brace,mirror}, line width=0.8pt] (-8.2,-16) -- (-3,-16) node[midway, below=30pt] {};
    \coordinate (start7) at (-5.6,-16.1);
    \coordinate (bend7) at (-5.6,-17);

    \draw[->, thin] (start7) -- (bend7)node[pos=0.8, above left, xshift=20pt] {}; 

\end{tikzpicture}
\end{center}

\section{Preliminaries}\label{pre}
We first present some basic notations and useful lemmas.
In this paper, the following notations
$$x_0=t,\ x=(x_1,x_2),\ \partial\in\{\partial_t,\partial_1,\partial_2\},\ \omega=\left(\frac{x_1}{r},\frac{x_2}{r}\right),\ \omega^\bot
=\left(\frac{x_2}{r},-\frac{x_1}{r}\right)$$
will be used.
Let $\Gamma$ be one of the vector fields
\begin{equation}
\partial_i,\ S=t\partial_t+r\partial_r, \ \Omega=x_2\partial_1-x_1\partial_2,
\end{equation}
where $i=0,1,2$.
Due to the communication relationship of $\Gamma$ and $\Box$, for multiple index $\alpha\in\mathbb{N}_+^5$, equation \eqref{phieq} is reduced to
\begin{equation}\label{Gphieq}
  \Box \Gamma^\alpha\phi=\sum_{|\beta+\gamma|\leq|\alpha|}C_\alpha^\beta g^{ijk}\partial_i\Gamma^\beta\phi \partial_j\partial_k\Gamma^\gamma\phi+\sum_{|\beta+\gamma+\delta|\leq|\alpha|}C_\alpha^\beta C_{\gamma+\delta}^\gamma h^{ijkl}\partial _i\Gamma^\beta\phi\partial_j\Gamma^\gamma\phi\partial_k\partial_l\Gamma^\delta\phi,
\end{equation}
with initial data
 \begin{equation}\label{Gphidata}
t=0:\quad \Gamma^\alpha\phi(0,x)=\epsilon \Gamma^\alpha\phi_0(x),\quad \partial_t\Gamma^\alpha\phi(0,x)=\epsilon \Gamma^\alpha\phi_1(x).
 \end{equation}


Based on these weighted vector fields and Sobolev embedding theorem, the following weighted $L^\infty-L^2$ decay estimates hold, see Sideris \cite{S}.
\begin{lem}\label{lem0}
For any $\phi\in C^\infty([0,T]\times\mathbb{R}^2)$, the following Sobolev type inequalities hold:
\begin{align}
      &\langle r\rangle^{1/2}| \phi(t,x)|\lesssim \sum_{|\alpha|\leq2}\|\Gamma^\alpha \phi\|_{L^2},\label{24}\\
      &\langle r\rangle^{1/2}\langle t-r\rangle^{1/2}|\partial\phi(t,x)|\lesssim \sum_{|\alpha|\leq1}\Big(\|\partial\Omega^\alpha \phi\|_{L^2}+\|\langle t-r\rangle\partial^2\Omega^\alpha \phi\|_{L^2}\Big),\\
      &\langle r\rangle^{1/2}\langle t-r\rangle|\partial^2 \phi(t,x)|\lesssim \sum_{|\alpha|\leq2} \|\langle t-r\rangle\partial^2\Gamma^\alpha \phi\|_{L^2},\label{4.3.0}
 \end{align}
 provided that the norms on the right-hand side are finite. 
\end{lem}

If the wave $\phi$ is compactly supported, the following Hardy type inequality helps us gain one derivative with respect to space or time. 
\begin{lem}\label{lemhardy}
Let $\phi\in C^\infty([0,T]\times\mathbb{R}^2)$. If $\phi$ is supported in $|x|\leq t+1$, then we have
  \begin{equation}\label{1.11}
    \|\langle t-r\rangle^{-1} \phi\|_{L^2}\lesssim \|\partial \phi\|_{L^2},
  \end{equation}
  and
  \begin{equation}\label{1.12}
    \|\langle r\rangle^{1/2}\langle t-r\rangle^{-1} \phi\|_{L^\infty}\lesssim \sum_{|\alpha|\leq2}\|\partial\Gamma^\alpha \phi\|_{L^2}.
  \end{equation}
\end{lem}
\emph{Proof.} For the Hardy type inequality \eqref{1.11}, see e.g. Lemma 3.3 in Katayama \cite{Kata} and Lemma 4.4 in Hidano \cite{Hida}. Combining \eqref{1.11} and Lemma \ref{lem0}, inequality \eqref{1.12} follows, see also Lemma 2.6 in Zha \cite{Zha}. \hfill$\Box$ 

When employing the null conditions outside the cone, i.e., $|x|\ge \langle t \rangle/2$, the following lemma is useful to control the $L^\infty$ norm of $\Gamma \phi$ by the energy.
\begin{lem}\label{lemout}
Let $\phi\in C^\infty([0,T]\times\mathbb{R}^2)$. If $\phi$ is supported in $|x|\leq t+1$, then we have
  \begin{equation}\label{Sobout}
    \left\|\Gamma^\alpha\phi\right\|_{L_x^\infty\left(|x|\ge \langle t \rangle/2\right)}\lesssim E_{\kappa}^{\frac 12}(\phi(t)),~~~|\alpha|+2\le \kappa.
  \end{equation}
\end{lem}
\emph{Proof.} This lemma comes from inequality (2.40) in Lemma 2.8 in Peng and Zha \cite{PZ}. \hfill$\Box$

We then introduce a $k$-th order Morawetz type energy
\begin{equation}\label{321}
\begin{split}
  \tilde{M}_\kappa(\phi)=& \sum_{|\alpha|\leq \kappa-1}\int_{\mathbb{R}^2}\frac{t+r+1}{8}\left(\partial_t\Gamma^\alpha \phi+\partial_r\Gamma^\alpha \phi\right)^2dx\\
  &+ \sum_{|\alpha|\leq \kappa-1}\int_{\mathbb{R}^2}\frac{t-r+1}{8}\left(\partial_t\Gamma^\alpha \phi-\partial_r\Gamma^\alpha \phi\right)^2dx\\
   &+ \sum_{|\alpha|\leq \kappa-1}\int_{\mathbb{R}^2}\frac{t+1}{2r^2}\left(\Omega \Gamma^\alpha \phi\right)^2dx+ \sum_{|\alpha|\leq \kappa-1}\int_{\mathbb{R}^2}\frac{1}{8(t+1)}\left|\Gamma^\alpha \phi\right|^2dx,
\end{split}
\end{equation}
the construction of which is motivated by the Morawetz type weight \eqref{multi} below, and the details will be discussed in Section \ref{Moen}. We have the following  weighted $L^2-L^2$ estimates.
\begin{lem}\label{lemM1}
For any $\phi\in C^\infty([0,T]\times\mathbb{R}^2)$ supported in $\{x: |x|\leq t+1\}$, it holds
  \begin{equation}\label{210}
    \|\langle t-r\rangle^{1/2}\partial  \phi\|_{L^2}\lesssim \tilde{M}_1^{1/2}( \phi).
  \end{equation}
\end{lem}
\emph{Proof.} On one hand, noting that
\begin{equation}
  \nabla =\omega\partial_r+\frac{\omega^\bot}{r}\Omega,
\end{equation}
we have
\begin{equation}
  \begin{split}
    \int_{r<t+1}\langle t-r\rangle|\nabla  \phi|^2dx\lesssim& \int_{r<t+1}(t-r+1)|\nabla  \phi|^2dx\\
    \lesssim&\int_{r<t+1}(t-r+1)|\partial_r  \phi|^2dx+\int_{r<t+1}\frac{t-r+1}{r^2}|\Omega  \phi|^2dx.
  \end{split}
\end{equation}
Since $t-r+1\leq t+1$, we then obtain
\begin{equation}
  \begin{split}
    \int_{r<t+1}\langle t-r\rangle|\nabla  \phi|^2dx\lesssim&\int_{\mathbb{R}^2}(t-r+1)|\partial_t\phi+\partial_r\phi-\partial_t\phi+\partial_r\phi |^2dx+\int_{\mathbb{R}^2}\frac{t+1}{r^2}|\Omega  \phi|^2dx\\
    \lesssim&\int_{\mathbb{R}^2}(t-r+1)|\partial_t\phi+\partial_r\phi|^2+\int_{\mathbb{R}^2}(t-r+1)|\partial_t\phi-\partial_r\phi|^2dx\\
    &+\int_{\mathbb{R}^2}\frac{t+1}{r^2}|\Omega  \phi|^2dx\\
    \lesssim&\int_{\mathbb{R}^2}(t+r+1)|\partial_t\phi+\partial_r\phi|^2+\int_{\mathbb{R}^2}(t-r+1)|\partial_t\phi-\partial_r\phi|^2dx\\
    &+\int_{\mathbb{R}^2}\frac{t+1}{r^2}|\Omega  \phi|^2dx\\
    \lesssim&\tilde{M}_1(\phi).
  \end{split}
\end{equation}

On the other hand, we have
\begin{equation}
  \begin{split}
    \int_{r<t+1}\langle t-r\rangle|\partial_t \phi|^2dx\lesssim&\int_{\mathbb{R}^2}(t-r+1)|\partial_t\phi+\partial_r\phi+\partial_t\phi-\partial_r\phi |^2dx\\
        \lesssim&\int_{\mathbb{R}^2}(t-r+1)|\partial_t\phi+\partial_r\phi|^2+\int_{\mathbb{R}^2}(t-r+1)|\partial_t\phi-\partial_r\phi|^2dx\\
    \lesssim&\int_{\mathbb{R}^2}(t+r+1)|\partial_t\phi+\partial_r\phi|^2+\int_{\mathbb{R}^2}(t-r+1)|\partial_t\phi-\partial_r\phi|^2dx\\
    \lesssim&\tilde{M}_1(\phi).
  \end{split}
\end{equation}
Hence inequality \eqref{210} follows. \hfill$\Box$


\section{Morawetz type energy}\label{Moen}

In this section, we establish a Morawetz type energy estimate for the Cauchy problem of the scalar wave equation 
\begin{empheq}[left=\empheqlbrace]{align}
&\Box \Phi=F,\label{3.1}\\
&\Phi(0,x)=\Phi_0(x),\ \partial_t\Phi(0,x)=\Phi_1(x)
\end{empheq}
in $\R^2$, by using a modified scaling Morawetz type multiplier
\begin{equation}\label{multi}
Z\Phi=\partial_t\Phi+(t+1)\partial_t\Phi+r\partial_r\Phi+\frac12\Phi.
\end{equation}
The key to prove Theorem \ref{thm1} by our method is to show the non-negativity of this Morawetz type energy with compactly supported initial data, which is also true for linear wave systems with multiple propagation speeds, and this would in turn to prove Theorem \ref{thmhar}. 
\begin{remark}
A similar multiplier was used in \cite{LaiZhou} for inhomogeneous wave equations in $\R^n (n\ge 4)$ without compact support assumption, while in \cite{Lai20, LaZ21} a similar multiplier was introduced for the wave equations in $\R^3$ but with compactly supported initial data. Very recently, the standard scaling Morawetz multiplier was used in \cite{LRX25} to establish the lower bound of lifespan for the quasilinear wave equation in exterior domain in $\R^2$.
\end{remark}

Multiplying the both sides of  \eqref{3.1} by \eqref{multi} and integrating on the whole space $\mathbb{R}^2$, we get
\begin{equation}
  \begin{split}
    &\frac12\frac{d}{dt}\int_{\mathbb{R}^2}\left(|\partial_t\Phi|^2+|\nabla \Phi|^2\right)dx\\
    &+\frac12\frac{d}{dt}\int_{\mathbb{R}^2}\Big[(t+1)|\partial_t\Phi|^2+(t+1)|\nabla \Phi|^2+2r\partial_r\Phi\partial_t\Phi+\Phi\partial_t\Phi\Big]dx\\
    =&\int_{\mathbb{R}^2}F\cdot Z\Phi dx.
  \end{split}
\end{equation}
Let 
\begin{equation}
E_1(\Phi)=\frac12\int_{\mathbb{R}^2}\left(|\partial_t\Phi|^2+c_i^2|\nabla \Phi|^2\right)dx
\end{equation}
and
\begin{equation}\label{morawetz}
M_1(\Phi)=\frac12\int_{\mathbb{R}^2}\Big[(t+2)|\partial_t\Phi|^2+(t+2)|\nabla \Phi|^2+2r\partial_r\Phi\partial_t\Phi+\Phi\partial_t\Phi\Big]dx,
\end{equation}
which is called the first-order Morawetz type energy.

We first show that the Morawetz type energy $M_1(\Phi)$ is non-negative if the initial data $(\Phi_0,\Phi_1)$ of \eqref{3.1} have compact support. Without loss of generality, we assume that
\begin{equation} \label{compact}
\textup{supp}\ \Phi_0=\textup{supp}\ \Phi_1=\{x:|x|\leq 1^-\}.
\end{equation} 
\begin{lem}\label{th1}
Suppose that $\Phi_0$ and $\Phi_1$ satisfy \eqref{compact}. Then the Morawetz type energy defined by \eqref{morawetz} is non-negative
\begin{equation}\label{positive}
  M_1(\Phi)\geq 0.
\end{equation}
In fact,
\begin{equation}\label{4.6}
\begin{split}
   M_1(\Phi)\geq &E_1(\Phi)+ \int_{\mathbb{R}^2}\frac{t+r+1}{8}(\partial_t\Phi+\partial_r\Phi)^2dx\\
   &+\int_{\mathbb{R}^2}\frac{t-r+1}{8}(\partial_t\Phi-\partial_r\Phi)^2dx\\
   &+\int_{\mathbb{R}^2}\frac{t+1}{2r^2}(\Omega \Phi)^2dx+\int_{\mathbb{R}^2}\frac{1}{8(t+1)}|\Phi|^2dx.
\end{split}
\end{equation}
\end{lem}

\emph{Proof.} Multiplying $\partial_t\Phi$ on the left side of \eqref{3.1}, integrating on $\mathbb{R}^2$,  we get 
 \begin{equation}
\int_{\mathbb{R}^2}\Box\Phi\partial_t\Phi dx=\frac{d}{dt}E_1(\Phi).
 \end{equation}
 
In the following, we show $M_1(\Phi)-E_1(\Phi)$ is non-negative. 
Since
\begin{equation}
  \nabla \Phi=\omega\partial_r\Phi+\frac{\omega^\bot}{r}\Omega \Phi,
\end{equation}
on one hand, we have 
\begin{equation}\label{310}
\begin{split}
 \frac{1}{2}(t+1)|\nabla \Phi|^2=&\frac{t+1}{2}\big|\omega\partial_r\Phi+\frac{\omega^\bot}{r}\Omega \Phi\big|^2\\
=&\frac{t+1}{2}|\partial_r\Phi|^2+\frac{t+1}{2r^2}|\Omega \Phi|^2.
  \end{split}
 \end{equation}
On the other hand, direct calculation implies
\begin{equation}\label{311}
  \begin{split}
&\frac{t+1}{2}|\partial_t\Phi|^2+\frac{t+1}{2}|\partial_r \Phi|^2+r\partial_r\Phi\partial_t\Phi\\
    =&\frac{t+r+1}{4}(\partial_t\Phi+\partial_r \Phi)^2+\frac{t-r+1}{4}(\partial_t\Phi-\partial_r\Phi)^2.
  \end{split}
\end{equation}
Combing \eqref{310} and \eqref{311} yields
\begin{equation}\label{04.15}
  \begin{split}
    &M_1(\Phi)-E_1(\Phi)\\
    =&\int_{\mathbb{R}^2}\frac{t+r+1}{4}(\partial_t\Phi+\partial_r \Phi)^2dx+\int_{\mathbb{R}^2}\frac{t-r+1}{4}(\partial_t\Phi-\partial_r \Phi)^2dx\\
    &+\int_{\mathbb{R}^2}\frac{t+1}{2r^2}|\Omega \Phi|^2dx+\int_{\mathbb{R}^2}\frac12\Phi\partial_t\Phi dx.
  \end{split}
\end{equation}
It is obvious that the first two terms on the right hand side of \eqref{04.15} are non-negative for $\text{supp}{\Phi}=\{x:|x|\leq t+1\}$. Next we estimate the last term of $M_1(\Phi)$.
Noting that
\begin{equation}
  \begin{split}
   ( t+1)\partial_t \Phi+r\partial_r\Phi=\frac{t+r+1}{2}(\partial_t\Phi+\partial_r \Phi)+\frac{t-r+1}{2}(\partial_t\Phi-\partial_r \Phi),
  \end{split}
\end{equation}
we have
\begin{equation}\label{04.17}
  \begin{split}
&\int_{\mathbb{R}^2}\frac12\Phi\partial_t\Phi dx\\  =&\int_{\mathbb{R}^2}\Phi\left[\frac{t+r+1}{4(t+1)}(\partial_t\Phi+\partial_r \Phi)+\frac{t-r+1}{4(t+1)}(\partial_t\Phi-\partial_r \Phi)\right]dx\\
    &-\int_{\mathbb{R}^2}\frac{1}{2(t+1)}\Phi r\partial_r\Phi dx\\
    =&\underbrace{\int_{\mathbb{R}^2}\frac{t+r+1}{4(t+1)}\Phi(\partial_t\Phi+\partial_r \Phi)dx}_{I_1}+\underbrace{\int_{\mathbb{R}^2}\frac{t-r+1}{4(t+1)}\Phi(\partial_t\Phi-\partial_r \Phi)dx}_{I_2}\\
    &+\int_{\mathbb{R}^2}\frac{1}{2(t+1)}\Phi^2dx,
  \end{split}
\end{equation}
where integration by parts is used in the last line above. In the following, we show that $I_1,I_2$ can be absorbed by the square terms in \eqref{04.15} and \eqref{04.17}.
For $I_1$, since $\text{supp} \Phi=\{|x|:|x|\leq t+1^-\}$, then for $|x|< t+1$ we have
\begin{equation}
  \begin{split}
    \sqrt{\frac{t+r+1}{8(t+1)}}\frac{4(t+1)}{t+r+1}=\sqrt{\frac{2(t+1)}{t+r+1}}\geq1,
  \end{split}
\end{equation}
which means that
\begin{equation}
  \begin{split}
    \frac{t+r+1}{4(t+1)}\leq \sqrt{\frac{t+r+1}{8(t+1)}}=2\sqrt{\frac{t+r+1}{8}}\sqrt{\frac{1}{4(t+1)}}.
  \end{split}
\end{equation}
Thus by Cauchy-Schwartz inequality, we get
\begin{equation}\label{04.20}
  \begin{split}
    &I_1=\int_{\mathbb{R}^2}\frac{t+r+1}{4(t+1)}\Phi(\partial_t\Phi+\partial_r \Phi)dx\\
    \geq&-\int_{\mathbb{R}^2}\frac{t+r+1}{8}(\partial_t\Phi+\partial_r \Phi)^2dx-\int_{\mathbb{R}^2}\frac{1}{4(t+1)}\Phi^2dx.
  \end{split}
\end{equation}
For $I_2$, similar analysis as above yields
\begin{equation}
  \begin{split}
    \sqrt{\frac{t-r+1}{16(t+1)}}\frac{4(t+1)}{t-r+1}=\sqrt{\frac{t+1}{t-r+1}}\geq\sqrt{\frac{t-r+1}{t-r+1}}=1,
  \end{split}
\end{equation}
which implies
\begin{equation}
  \begin{split}
    \frac{t-r+1}{4(t+1)}\leq \sqrt{\frac{t-r+1}{16(t+1)}}=2\sqrt{\frac{t-r+1}{8}}\sqrt{\frac{1}{8(t+1)}}.
  \end{split}
\end{equation}
Hence, we have
\begin{equation}\label{04.23}
  \begin{split}
    &I_2=\int_{\mathbb{R}^2}\frac{t-r+1}{4(t+1)}\Phi(\partial_t\Phi-\partial_r \Phi)dx\\
    \geq&-\int_{\mathbb{R}^2}\frac{t-r+1}{8}(\partial_t\Phi-\partial_r \Phi)^2dx-\int_{\mathbb{R}^2}\frac{1}{8(t+1)}\Phi^2dx.
  \end{split}
\end{equation}
Finally, inserting \eqref{04.17}, \eqref{04.20} and \eqref{04.23} into \eqref{04.15}, we get \eqref{4.6}. This completes the proof.\hfill$\Box$

For notation convenience, we denote
\begin{equation}\label{M1bar}
\begin{split}
  \tilde{M}_1(\Phi)=& \int_{\mathbb{R}^2}\frac{t+r+1}{8}(\partial_t\Phi+\partial_r\Phi)^2dx+\int_{\mathbb{R}^2}\frac{t-r+1}{8}(\partial_t\Phi-\partial_r\Phi)^2dx\\
   &+\int_{\mathbb{R}^2}\frac{t+1}{2r^2}(\Omega \Phi)^2dx+\int_{\mathbb{R}^2}\frac{1}{8(t+1)}|\Phi|^2dx.
\end{split}
\end{equation}
In Lemma \ref{th1},  it shows that $M_1(\Phi)$ has a lower bound $\tilde{M}_1(\Phi)+E_1(\Phi)$. In fact, $\tilde{M}_1(\Phi)+E_1(\Phi)$ is also an upper bound of the Morawetz type energy $M_1(\Phi)$.
\begin{lem}\label{lem2}
Let $M_1(\Phi)$ be defined as in \eqref{morawetz}. Then the Morawetz type energy $M_1(\Phi)$ is equivalent to $\tilde{M}_1(\Phi)+E_1(\Phi)$, i.e.,
\begin{equation}\label{eqvMM}
M_1(\Phi)\simeq\tilde{M}_1(\Phi)+E_1(\Phi).
\end{equation}
\end{lem}
\emph{Proof.} By Lemma \ref{th1}, we only need to show $M_1(\Phi)$ is bounded by $\tilde{M}_1(\Phi)+E_1(\Phi)$. From \eqref{04.15} and \eqref{04.17} it holds
\begin{equation}\label{322}
  \begin{split}
    M_1(\Phi)=&\frac12\int_{\mathbb{R}^2}(|\partial_t\Phi|^2+|\nabla \Phi|^2)dx\\
    &+\frac12\int_{\mathbb{R}^2}\Big[(t+1)|\partial_t \Phi|^2+(t+1)|\nabla  \Phi|^2+2r\partial_r \Phi\partial_t \Phi+\Phi\partial_t \Phi\Big]dx\\
    =&E_1(\Phi)+\int_{\mathbb{R}^2}\frac{t+r+1}{4}(\partial_t\Phi+\partial_r \Phi)^2dx\\
    &+\int_{\mathbb{R}^2}\frac{t-r+1}{4}(\partial_t\Phi-\partial_r \Phi)^2dx\\
    &+\int_{\mathbb{R}^2}\frac{t+1}{2r^2}|\Omega \Phi|^2dx+\int_{\mathbb{R}^2}\frac{1}{2(t+1)}(\Phi)^2dx\\
      &+\int_{\mathbb{R}^2}\Big[\underbrace{\frac{t+r+1}{4(t+1)}\Phi(\partial_t\Phi+\partial_r \Phi)+\frac{t-r+1}{4(t+1)}\Phi(\partial_t\Phi-\partial_r \Phi)}_{J}\Big]dx.
  \end{split}
\end{equation}
We then decompose $J$ to obtain
\begin{equation}\label{J}
  \begin{split}
 J= &\frac{\Phi}{2( t+1)^{1/2}}\cdot\Big[\frac{t+r+1}{2(t+1)^{1/2}}(\partial_t\Phi+\partial_r \Phi)\Big]+\frac{\Phi}{2( t+1)^{1/2}}\cdot\Big[\frac{t-r+1}{2(t+1)^{1/2}}(\partial_t\Phi-\partial_r \Phi)\Big]\\
 \leq&\frac{1}{4(t+1)}|\Phi|^2+\frac{(t+r+1)^2}{4(t+1)}(\partial_t\Phi+\partial_r \Phi)^2+\frac{1}{4( t+1)}|\Phi|^2+\frac{(t-r+1)^2}{4(t+1)}(\partial_t\Phi-\partial_r \Phi)^2.\\
  \end{split}
\end{equation}
Noting that $0<r\leq t+1$, it holds 
\begin{equation}\label{324}
\frac{(t+r+1)^2}{4(t+1)}\leq\frac{2(t+1)(t+r+1)}{4(t+1)}=\frac{t+r+1}{2}
\end{equation}
and
\begin{equation}\label{325}
\frac{(t-r+1)^2}{4(t+1)}\leq\frac{t-r+1}{4}.
\end{equation}
By combining \eqref{J}, \eqref{324} and \eqref{325}, we come to
\begin{equation}
 J\leq\frac{}{2(t+1)}|\Phi|^2+\frac{t+r+1}{2}(\partial_t\Phi+\partial_r \Phi)^2+\frac{t-r+1}{4}(\partial_t\Phi-\partial_r \Phi)^2,
\end{equation}
together with \eqref{322}, this leads to 
\begin{equation}
  \begin{split}
    M_1(\Phi)\leq&E_1(\Phi)+\int_{\mathbb{R}^2}\frac{3(t+r+1)}{4}(\partial_t\Phi+\partial_r \Phi)^2dx\\
    &+\int_{\mathbb{R}^2}\frac{t-r+1}{2}(\partial_t\Phi-\partial_r \Phi)^2dx\\
    &+\int_{\mathbb{R}^2}\frac{t+1}{2r^2}|\Omega \Phi|^2dx+\int_{\mathbb{R}^2}\frac{1}{t+1}(\Phi)^2dx.
  \end{split}
\end{equation}
Thus, we get the upper bound of the Morawetz type energy
\begin{equation}
M_1(\Phi)\lesssim E_1(\Phi)+\tilde{M}_1(\Phi),
\end{equation}
and hence the proof is completed by combining Lemma \ref{th1}.\hfill$\Box$

\section{Null conditions}
For nonlinear terms of equation \eqref{phieq}, the null conditions \eqref{fstnull} and \eqref{secnull} can improve the decay. This analysis was contained in Lemma 6.6.4 of  \cite{Hormander} by H\"ormander with the full generalized vector fields. With absence of Lorentz boost operator, the following lemma can be referred to Zha \cite{Zha}.
\begin{lem}\label{nulllem}
If $g^{ijk}$ satisfy the first null condition \eqref{fstnull}, then for any $|x|\geq t/2$,
 \begin{equation}
|g^{ijk}\partial_k\phi\partial_{ij}^2\psi|\lesssim \langle t+r\rangle^{-1}\left(|\Gamma \phi||\partial^2\psi|+|\partial\phi||\partial\Gamma\psi|+\langle t-r\rangle|\partial\phi||\partial^2\psi|\right),
\end{equation}
\begin{equation}
  \begin{split}
|g^{ijk}\partial_i\phi\partial_j\psi\partial_k\varphi|\lesssim \langle t+r\rangle^{-1}\Big(&|\Gamma \phi||\partial\psi||\partial\varphi|+|\partial\phi||\Gamma\psi||\partial\varphi|+|\partial\phi||\partial\psi||\Gamma\varphi|\\
&+\langle t-r\rangle|\partial\phi||\partial\psi||\partial\varphi|\Big).
  \end{split}
\end{equation}
If $h^{ijkl}$ satisfy the second null condition \eqref{secnull}, then for any $|x|\geq t/2$,
\begin{equation}
  \begin{split}
|h^{ijkl}\partial_k\phi\partial_l\varphi\partial_{ij}^2\psi|\lesssim \langle t+r\rangle^{-1}\Big(&|\Gamma \phi||\partial\varphi||\partial^2\psi|+|\partial\phi||\Gamma\varphi||\partial^2\psi|+|\partial\phi||\partial\varphi||\partial\Gamma\psi|\\
&+\langle t-r\rangle|\partial\phi||\partial\varphi||\partial^2\psi|\Big),
  \end{split}
\end{equation}
\begin{equation}
  \begin{split}
|h^{ijkl}\partial_i\phi\partial_j\psi\partial_k\varphi\partial_l\chi|\lesssim \langle t+r\rangle^{-1}\Big(&|\Gamma \phi||\partial\psi||\partial\varphi||\partial\chi|+|\partial\phi||\Gamma\psi||\partial\varphi||\partial\chi|\\
&+|\partial\phi||\partial\psi||\Gamma\varphi||\partial\chi|+|\partial\phi||\partial\psi||\partial\varphi||\Gamma\chi|\\
&+\langle t-r\rangle|\partial\phi||\partial\psi||\partial\varphi||\partial\chi|\Big).
  \end{split}
\end{equation}
\end{lem}

If only the spatial rotational operator $\Omega$ is involved, then the null conditions offer $r^{-1}$ decay.
\begin{lem}\label{newnull}
If $g^{ijk}$ satisfy the first null condition \eqref{fstnull}, then 
\begin{equation}\label{45}
|g^{ijk}\partial_i\phi\partial_j\psi\partial_k\varphi|\lesssim r^{-1}\left(|\Omega \phi||\partial\psi||\partial\varphi|+|\partial\phi||\Omega\psi||\partial\varphi|+|\partial\phi||\partial\psi||\Omega\varphi|\right).
\end{equation}
If $h^{ijkl}$ satisfy the second null condition \eqref{secnull}, then 
\begin{equation}\label{46}
  \begin{split}
|h^{ijkl}\partial_i\phi\partial_j\psi\partial_k\varphi\partial_l\chi|\lesssim r^{-1}\Big(&|\Omega \phi||\partial\psi||\partial\varphi||\partial\chi|+|\partial\phi||\Omega\psi||\partial\varphi||\partial\chi|\\
&+|\partial\phi||\partial\psi||\Omega\varphi||\partial\chi|+|\partial\phi||\partial\psi||\partial\varphi||\Omega\chi|\Big).
  \end{split}
\end{equation}
\end{lem} 
\emph {Proof.} This lemma can be proved by direct calculation. We present the detailed proof for inequality \eqref{45}, and \eqref{46} can be shown analogously.
Noting that $$\partial_i=\omega_i\partial_r+\frac{\omega^\bot_i}{r}\Omega,$$ 
 we have
\begin{equation}
\begin{split}
|rg^{ijk}\partial_i\phi\partial_j\psi\partial_k\varphi|\leq &|rg^{ijk}\omega_i\partial_r\phi\partial_j\psi\partial_k\varphi|+|g^{ijk}\omega^\bot_i\Omega\phi\partial_j\psi\partial_k\varphi|\\
\leq&|rg^{ijk}\omega_i\omega_j\partial_r\phi\partial_r\psi\partial_k\varphi|+|g^{ijk}\omega_i\omega^\bot_j\partial_r\phi\Omega\psi\partial_k\varphi|\\
&+|g^{ijk}\omega^\bot_i\Omega\phi\partial_j\psi\partial_k\varphi|\\
\leq&|rg^{ijk}\omega_i\omega_j\omega_k\partial_r\phi\partial_r\psi\partial_r\varphi|+|g^{ijk}\omega_i\omega_j\omega^\bot_k\partial_r\phi\partial_r\psi\Omega\varphi|\\
&+|g^{ijk}\omega_i\omega^\bot_j\partial_r\phi\Omega\psi\partial_k\varphi|+|g^{ijk}\omega^\bot_i\Omega\phi\partial_j\psi\partial_k\varphi|.
\end{split}
\end{equation}
Since $g^{ijk}\omega_i\omega_j\omega_k=0$ and $\partial_r=\omega_1\partial_1+\omega_2\partial_2$, inequality \eqref{45} is follows. \hfill$\Box$

The following Klainerman-Sideris inequality makes up the absence of the Lorentz boost operator. 
\begin{lem}[Klainerman-Sideris \cite{KlainSid}]
Let $\phi\in C^\infty([0,T]\times\mathbb{R}^2)$. Then
  \begin{equation}\label{ksineq}
  \|\langle t-r\rangle\partial\nabla \phi\|_{L^2}\lesssim E_2^{1/2}(\phi)+\|(t+r)\Box\phi\|_{L^2}.
  \end{equation}
\end{lem}

Define
\begin{equation}
\mathcal{X}_\mu(\phi):=\sum_{|\alpha|\leq\mu-2} \|\langle t-r\rangle\partial\nabla \Gamma^\alpha\phi\|_{L^2},
\end{equation}
then we have
\begin{lem}\label{lem3.3}
Let $\mu\geq 10$, and let $\phi$ be the solution to \eqref{phieq}-\eqref{phidata}. If $E_\mu^{1/2}(\phi)=\varepsilon_0$ is sufficiently small for any $t\in[0,T]$, then for $t\in[0,T]$ it holds
\begin{equation}\label{3.11}
\mathcal{X}_\mu(\phi)\lesssim E_\mu^{1/2}(\phi).
\end{equation}
\end{lem}
\emph{Proof.} The proof was given in Sideris-Tu \cite{st}.

%

\section{Energy estimates}
 Let $\kappa\geq12$ and $\mu=\kappa-1$. In this section, we obtain two energy estimates with enough time decay for the Cauchy problem of \eqref{phieq}. In detail, we are going to establish:
\begin{equation}
  \begin{split}
        \frac{d}{dt}M_\kappa( \phi(t))\lesssim\langle t\rangle^{-1}E_{\mu}^{1/2}( \phi(t))M_\kappa( \phi(t)), \qquad 0<t<T,
  \end{split}
\end{equation}
and 
\begin{equation}
 \frac{d}{dt}E_{\mu}(\phi(t))\lesssim\langle t\rangle^{-3/2}E_\mu(\phi(t))M_\kappa^{1/2}(\phi(t)), \qquad 0<t<T.
\end{equation}
Then Theorem \ref{thm1} can be proved by a bootstrap argument getting $E_\mu^{1/2}(\phi(t))\lesssim\epsilon$ and $M_\kappa^{1/2}(\phi(t))\lesssim(1+t)^\sigma$ for any $t>0$ and some small constant $0<\sigma<1$.

\subsection{Morawetz energy estimate}
Let $|\alpha|\leq \kappa-1$. Multiplying $$Z\Gamma^\alpha\phi=\partial_t\Gamma^\alpha\phi+(t+1)\partial_t\Gamma^\alpha \phi+r\partial_r\Gamma^\alpha \phi+\frac12\Gamma^\alpha \phi$$ on both sides of equation \eqref{Gphieq}, then integrating them on $\mathbb{R}^2$, we get
\begin{equation}
  \begin{split}
    \frac{d}{dt}M_\kappa(\phi) = &\sum_{|\alpha|=0}^{\kappa-1}\underbrace{\int_{\mathbb{R}^2} (g^{ijk}\partial_i\phi \partial_j\partial_k\Gamma^\alpha\phi +h^{ijkl}\partial _i\phi \partial_j\phi\partial_k\partial_l\Gamma^\alpha\phi)Z\Gamma^\alpha\phi dx}_{L_1(\mathbb{R}^2)}\\
    &+\sum_{|\alpha|=0}^{\kappa-1}\underbrace{\sum_{|\beta+\gamma|\leq|\alpha|\atop |\beta|, |\gamma|<|\alpha|}C_\alpha^\beta\int_{\mathbb{R}^2} g^{ijk}\partial_i\Gamma^\beta\phi \partial_j\partial_k\Gamma^\gamma\phi Z\Gamma^\alpha\phi dx}_{L_2(\mathbb{R}^2)}\\
    &+\sum_{|\alpha|=0}^{\kappa-1}\underbrace{\sum_{|\beta+\gamma+\delta|\leq|\alpha|\atop |\beta|, |\gamma|, |\delta|<|\alpha|}C_\alpha^\beta C_{\gamma+\delta}^\gamma\int_{\mathbb{R}^2} h^{ijkl}\partial _i\Gamma^\beta\phi\partial_j\Gamma^\gamma\phi\partial_k\partial_l\Gamma^\delta\phi Z\Gamma^\alpha\phi dx}_{L_3(\mathbb{R}^2)}.
  \end{split}
\end{equation}
By using integration by parts, we get
\begin{equation}
  \begin{split}
L_1(\mathbb{R}^2)=&-\frac12\frac{d}{dt}\int_{\mathbb{R}^2} (t+2)\left(g^{ijk}\partial_i\phi \partial_{j}\Gamma^\alpha\phi \partial_k\Gamma^\alpha\phi+h^{ijkl}\partial _i\phi\partial_j\phi\partial_{k}\Gamma^\alpha\phi \partial_l\Gamma^\alpha\phi\right) dx\\
&\underbrace{+\frac12\int_{\mathbb{R}^2} g^{ijk}(S+2)\partial_i\phi \partial_{j}\Gamma^\alpha\phi \partial_k\Gamma^\alpha\phi +h^{ijkl}(S+2)(\partial _i\phi\partial_j\phi)\partial_{k}\Gamma^\alpha\phi \partial_l\Gamma^\alpha\phi dx}_{L_{11}(\mathbb{R}^2)}\\
&\underbrace{-\int_{\mathbb{R}^2} [g^{ijk}\partial _i\partial_k\phi\partial_j\Gamma^\alpha\phi - h^{ijkl}\partial _l(\partial_i\phi\partial_j\phi)\partial_k\Gamma^\alpha\phi] Z\Gamma^\alpha\phi dx}_{L_{12}(\mathbb{R}^2)}.
  \end{split}
\end{equation}
Let
\begin{equation}
  \begin{split}
    \bar{M}_\kappa( \phi)=M_\kappa( \phi)+\sum_{|\alpha|\leq\kappa-1}\frac12\int_{\mathbb{R}^2} (t+2)\left(g^{ijk}\partial_i\phi \partial_{j}\Gamma^\alpha\phi \partial_k\Gamma^\alpha\phi+h^{ijkl}\partial _i\phi\partial_j\phi\partial_{k}\Gamma^\alpha\phi \partial_l\Gamma^\alpha\phi\right) dx,
  \end{split}
\end{equation}
then it holds
\begin{equation}\label{4.5}
  \begin{split}
    \frac{d}{dt} \bar{M}_\kappa(\phi)\lesssim L_{11}+L_2+L_3+L_{12}.
  \end{split}
\end{equation}

We first show $\bar{M}_\kappa(\phi)$ is equivalent to $M_\kappa(\phi)$. Divide $\mathbb{R}^2$ into two regions: $D_1=\{x: |x|\leq\frac{t}{2}\ \text {or}\ 2t\leq|x|\leq t+1\}$ and $D_2=\{x:\frac{t}2<|x|<2t\}$. For $D_1$, based on $|x|\leq t+1$, there is a uniform constant $C$ such that $$t+2\leq C\langle t-r\rangle,$$
which result in
\begin{equation}\label{wave57}
   \begin{split}
  & \int_{D_1} (t+2)\left(g^{ijk}\partial_k\phi +h^{ijkl}\partial _k\phi\partial_l\phi\right)\partial_{j}\Gamma^\alpha\phi \partial_i\Gamma^\alpha\phi dx\\
   \lesssim&\left\|\partial\phi+(\partial\phi)^2\right\|_{L^\infty}\left\|\langle t-r\rangle^{1/2}\partial\Gamma^\alpha \phi\right\|_{L^2}\left\|\langle t-r\rangle^{1/2}\partial\Gamma^\alpha \phi\right\|_{L^2}.
  \end{split}
\end{equation}
By combining Lemma \ref{lem0} and Lemma \ref{lemM1}, we  further have
\begin{equation}\label{56}
   \begin{split}
   &\int_{D_1} (t+2)\left(g^{ijk}\partial_k\phi +h^{ijkl}\partial _k\phi\partial_l\phi\right)\partial_{j}\Gamma^\alpha\phi \partial_i\Gamma^\alpha\phi dx\\
   \lesssim &\left[E_3^{1/2}(\phi)+E_3(\phi)\right]\tilde{M}_\kappa(\phi)
   \lesssim\varepsilon \tilde{M}_\kappa(\phi),
  \end{split}
\end{equation}
provided $E_3^{1/2}\leq \epsilon\ll1$. For $D_2$, it holds that $r\sim t$. We then  employ  Lemma \ref{newnull}, the Hardy type inequality \eqref{1.12}, inequalities \eqref{24} and \eqref{Sobout} to get 
\begin{equation}\label{57}
   \begin{split}
  & \int_{D_2} (t+2)\left(g^{ijk}\partial_k\phi +h^{ijkl}\partial _k\phi\partial_l\phi\right)\partial_{j}\Gamma^\alpha\phi \partial_i\Gamma^\alpha\phi dx\\
   \lesssim&\left\|\partial\phi+(\partial\phi)^2\right\|_{L^\infty}\left\|\partial\Gamma^\alpha \phi\right\|_{L^2}^2+\left\|(\Omega\phi+\partial\phi\Omega\phi)(\partial\Gamma^\alpha \phi)^2\right\|_{L^1}\\
   &+\left\|\left[\partial\phi+(\partial\phi)^2\right]\partial\Gamma^\alpha \phi\Omega\Gamma^\alpha \phi\right\|_{L^1(D_2)}\\
   \lesssim&\left[E_3^{1/2}(\phi)+E_3(\phi)\right]E_\kappa(\phi)+\left\|\Omega\phi+\partial\phi\Omega\phi\right\|_{L^\infty}\left\|\partial\Gamma^\alpha \phi\right\|_{L^2}^2\\
   &+\left\|r^{1/2}\left[\partial\phi+(\partial\phi)^2\right]\right\|_{L^\infty}\left\|\partial\Gamma^\alpha \phi\right\|_{L^2}\left\|\frac{1}{r^{1/2}}\Omega\Gamma^\alpha\phi\right\|_{L^2(D_2)}\\
   \lesssim&\varepsilon E_\kappa(\phi)+\varepsilon E_\kappa^{1/2}(\phi)\tilde{M}_\kappa^{1/2}(\phi),
\end{split}
\end{equation}
 provided $E_3^{1/2}\leq \epsilon\ll1$.
Combing \eqref{56} and \eqref{57}, we get
\begin{equation}
  \bar{M}_\kappa( \phi)\thickapprox M_\kappa( \phi).
\end{equation}

Now we estimate each term on the right hand side of \eqref{4.5}. As discussed above, we consider them in two cases: $D_1=\{x: |x|\leq\frac{t}{2}\ \text {or}\ 2t\leq|x|\leq t+1\}$ and $D_2=\{x:\frac{t}2<|x|<2t\}$. 

$\bullet$ \emph{Case \romannumeral1:} In region $D_1=\{x: |x|\leq\frac{t}{2}\ \text {or}\ 2t\leq|x|\leq t+1\}$. \\
In this case, we use the fact that 
\begin{equation}
1\lesssim\frac{\langle t-r\rangle}{\langle t\rangle}.
\end{equation}
For  $L_{11}(D_1)$, since $t\partial_t+r\partial_r=S$ is one of the $\Gamma$ operators, and by Lemma \ref{lem0} and Lemma \ref{lemM1} it is easy to get
\begin{equation}\label{510}
\begin{split}
    L_{11}(D_1)\lesssim&\langle t\rangle^{-1} \left\|\partial\Gamma\phi+\partial\phi\partial\Gamma\phi\right\|_{L^\infty}\left\|\langle t-r\rangle^{1/2}\partial\Gamma^\alpha \phi\right\|_{L^2}^2\\
    \lesssim& \langle t\rangle^{-1}E_\mu^{1/2}(\phi)M_\kappa(\phi).
\end{split}
\end{equation}
For $L_2(D_1)$, if $|\beta|\leq |\alpha|/2$, 
noting that
\begin{equation}\label{StoM}
t\partial_t+r\partial_r=\frac12(t+r)(\partial_t+\partial_r)+\frac12(t-r)(\partial_t-\partial_r)
\end{equation}
and 
 \begin{equation}\label{t+r}
 \langle t+r\rangle^{1/2}\lesssim \langle r\rangle^{1/2}\langle t-r\rangle^{1/2},
 \end{equation}
we then have 
\begin{equation}
\begin{split}
L_2(D_1)\lesssim&\langle t\rangle^{-1}\left\|\partial\Gamma^\beta\phi\right\|_{L^\infty}\left\|\langle t-r\rangle\partial^2\Gamma^\gamma\phi\right\|_{L^2}\left\|\partial\Gamma^\alpha\phi\right\|_{L^2}\\
&+\langle t\rangle^{-1}\left\|\langle t-r\rangle \partial \Gamma^\beta\phi\partial^2\Gamma^\gamma\phi(t-r)(\partial_t-\partial_r)\Gamma^\alpha\phi\right\|_{L^1}\\
&+\langle t\rangle^{-1}\left\|\langle t-r\rangle \partial \Gamma^\beta\phi\partial^2\Gamma^\gamma\phi(t+r)(\partial_t+\partial_r)\Gamma^\alpha\phi\right\|_{L^1}\\
&+\langle t\rangle^{-1}\left\|\langle t-r\rangle \partial \Gamma^\beta\phi\partial^2\Gamma^\gamma\phi\Gamma^\alpha\phi\right\|_{L^1}\\
\lesssim&\langle t\rangle^{-1}\left\|\partial\Gamma^\beta\phi\right\|_{L^\infty}\left\|\langle t-r\rangle\partial^2\Gamma^\gamma\phi\right\|_{L^2}\left\|\partial\Gamma^\alpha\phi\right\|_{L^2}\\
&+\langle t\rangle^{-1}\left\|\langle t-r\rangle^{1/2}\partial \Gamma^\beta\phi\right\|_{L^\infty}\left\|\langle t-r\rangle\partial^2\Gamma^\gamma\phi\right\|_{L^2}\left\|\langle t-r\rangle^{1/2}(\partial_t-\partial_r)\Gamma^\alpha\phi\right\|_{L^2}\\
&+\langle t\rangle^{-1}\left\|\langle r\rangle^{1/2}\langle t-r\rangle^{1/2} \partial \Gamma^\beta\phi\right\|_{L^\infty}\left\|\langle t-r\rangle\partial^2\Gamma^\gamma\phi\right\|_{L^2}\left\|\langle t+r\rangle^{1/2}(\partial_t+\partial_r)\Gamma^\alpha\phi\right\|_{L^2}\\
&+\langle t\rangle^{-1}\left\|\langle t-r\rangle^{1/2} \partial \Gamma^\beta\phi\right\|_{L^\infty}\left\|\langle t-r\rangle\partial^2\Gamma^\gamma\phi\right\|_{L^2}\left\|\frac{1}{\langle t\rangle^{1/2}}\Gamma^\alpha\phi\right\|_{L^2}\\
\lesssim& \langle t\rangle^{-1} E_\mu^{1/2}(\phi)\mathcal{X}_\kappa(\phi)E_\kappa^{1/2}(\phi)+\langle t\rangle^{-1}E_\mu^{1/2}(\phi)\mathcal{X}_\kappa(\phi)M_\kappa^{1/2}(\phi)\\
\lesssim&\langle t\rangle^{-1} E_\mu^{1/2}(\phi)E_\kappa^{1/2}(\phi)\left[E_\kappa^{1/2}(\phi)+M_\kappa^{1/2}(\phi)\right],
\end{split}
\end{equation}
while if $|\gamma|\leq |\alpha|/2$, we have
\begin{equation}
\begin{split}
L_2(D_1)\lesssim&\langle t\rangle^{-1}\left\|\partial\Gamma^\beta\phi\right\|_{L^2}\left\|\langle t-r\rangle\partial^2\Gamma^\gamma\phi\right\|_{L^\infty}\left\|\partial\Gamma^\alpha\phi\right\|_{L^2}\\
&+\langle t\rangle^{-1}\left\|\langle t-r\rangle \partial \Gamma^\beta\phi\partial^2\Gamma^\gamma\phi(t-r)(\partial_t-\partial_r)\Gamma^\alpha\phi\right\|_{L^1}\\
&+\langle t\rangle^{-1}\left\|\langle t-r\rangle \partial \Gamma^\beta\phi\partial^2\Gamma^\gamma\phi(t+r)(\partial_t+\partial_r)\Gamma^\alpha\phi\right\|_{L^1}\\
&+\langle t\rangle^{-1}\left\|\langle t-r\rangle \partial \Gamma^\beta\phi\partial^2\Gamma^\gamma\phi\Gamma^\alpha\phi\right\|_{L^1}\\
\lesssim&\langle t\rangle^{-1}\left\|\partial\Gamma^\beta\phi\|_{L^2}\|\langle t-r\rangle\partial^2\Gamma^\gamma\phi\right\|_{L^\infty}\left\|\partial\Gamma^\alpha\phi\right\|_{L^2}\\
&+\langle t\rangle^{-1}\left\|\langle t-r\rangle^{1/2}\partial \Gamma^\beta\phi\right\|_{L^2}\left\|\langle t-r\rangle\partial^2\Gamma^\gamma\phi\right\|_{L^\infty}\left\|\langle t-r\rangle^{1/2}(\partial_t-\partial_r)\Gamma^\alpha\phi\right\|_{L^2}\\
&+\langle t\rangle^{-1}\left\|\langle t-r\rangle^{1/2} \partial \Gamma^\beta\phi\right\|_{L^2}\left\|\langle r\rangle^{1/2}\langle t-r\rangle\partial^2\Gamma^\gamma\phi\right\|_{L^\infty}\left\|\langle t+r\rangle^{1/2}(\partial_t+\partial_r)\Gamma^\alpha\phi\right\|_{L^2}\\
&+\langle t\rangle^{-1}\left\|\langle t-r\rangle^{1/2} \partial \Gamma^\beta\phi\right\|_{L^2}\left\|\langle t-r\rangle\partial^2\Gamma^\gamma\phi\right\|_{L^\infty}\left\|\frac{1}{\langle t\rangle^{1/2}}\Gamma^\alpha\phi\right\|_{L^2}\\
\lesssim& \langle t\rangle^{-1} E_\kappa^{1/2}(\phi)\mathcal{X}_\mu(\phi)E_\kappa^{1/2}(\phi)+\langle t\rangle^{-1}M_\kappa^{1/2}(\phi)\left[E_\mu^{1/2}(\phi)+\mathcal{X}_\mu(\phi)\right]M_\kappa^{1/2}(\phi)\\
\lesssim&\langle t\rangle^{-1} E_\mu^{1/2}(\phi)\left[E_\kappa(\phi)+M_\kappa(\phi)\right].
\end{split}
\end{equation}
We then conclude that
\begin{equation}\label{517}
\begin{split}
L_2(D_1)\lesssim \langle t\rangle^{-1} E_\mu^{1/2}(\phi)\left[E_\kappa(\phi)+M_\kappa(\phi)\right].
\end{split}
\end{equation}
At last, for $L_3$, provided $E_{[|\alpha|/2]}\ll 1$, we may divide the estimate into two cases: $|\delta|\geq|\alpha|/2$ and $|\delta|\leq |\alpha|/2$, and similar discuss as above yields
\begin{equation}
\begin{split}
L_3(D_1)\lesssim \langle t\rangle^{-1} E_\mu^{1/2}(\phi)\left[E_\kappa(\phi)+M_\kappa(\phi)\right].
\end{split}
\end{equation}
In summary, in the region of $D_1$, we obtain
\begin{equation}\label{518}
\begin{split}
  L_{11}(D_1)+ L_{2}(D_1)+L_{3}(D_1)
  \lesssim \langle t\rangle^{-1} E_\mu^{1/2}(\phi)\left[E_\kappa(\phi)+M_\kappa(\phi)\right].
  \end{split}
\end{equation}

$\bullet$ \emph{Case \romannumeral2:} In region $D_2=\{x:\frac{t}2<|x|<2t\}$.\\
In this case, it holds that
\begin{equation}
\langle t\rangle\sim \langle r\rangle.
\end{equation}
To avoid redundant computation, here we only give the detailed calculation for the quadratic terms with coefficients $g^{ijk}$. The cubic terms with coefficients $h^{ijkl}$ can be estimated similarly, except for one more $L^\infty$ norm of the solution with lower-order derivatives. And by Sobolev inequality, this $L^\infty$ norm is bounded by $E_{[|\alpha|/2]}$, which is much smaller than one.
Hence, by Lemma \ref{newnull}, we have
\begin{equation}\label{520}
\begin{split}
  &L_{11}(D_2)+L_{2}(D_2)+ L_{3}(D_2)\\
  \lesssim &\langle t\rangle^{-1} \left(\underbrace{\|\Omega\Gamma\phi(\partial\Gamma^\alpha\phi)^2\|_{L^1(D_2)}}_{l_1}+\underbrace{\|\partial\Gamma\phi\partial\Gamma^\alpha\phi\Omega\Gamma^\alpha\phi\|_{L^1(D_2)}}_{l_2}\right)\\
  &+\langle t\rangle^{-1} \left(\underbrace{\|\Omega\Gamma^\beta\phi\partial^2\Gamma^\gamma\phi Z\Gamma^\alpha\phi\|_{L^1(D_2)}}_{l_3}+\underbrace{\|\partial\Gamma^\beta\phi\partial\Omega\Gamma^\gamma\phi Z\Gamma^\alpha\phi\|_{L^1(D_2)}}_{l_4}\right).
  \end{split}
\end{equation}
For $l_1$, by employing the Hardy type inequality, it holds
\begin{equation}\label{521}
\begin{split}
l_1\lesssim \left\|\langle t-r\rangle^{-1}\Omega\Gamma\phi\right\|_{L^\infty}\left\|\langle t-r\rangle^{1/2}\partial\Gamma^\alpha\phi\right\|_{L^2}^2\lesssim E_\mu^{1/2}(\phi)M_\kappa(\phi).
\end{split}
\end{equation}
For $l_2$, we estimate it as
\begin{equation}\label{522}
l_2\lesssim \left\|\langle r\rangle^{1/2}\partial\Gamma\phi\right\|_{L^\infty}\left\|\partial\Gamma^\alpha\phi\right\|_{L^2}\left\|\frac{1}{\langle r\rangle^{1/2}}\Omega\Gamma^\alpha\phi\right\|_{L^2}\lesssim E_\mu^{1/2}(\phi) E_\kappa^{1/2}(\phi)M_\kappa^{1/2}(\phi).
\end{equation}
For $l_3$ and $l_4$, by \eqref{StoM}, we have
\begin{equation}\label{523}
\begin{split}
l_3\lesssim& \underbrace{\left\|\Omega\Gamma^\beta\phi\partial^2\Gamma^\gamma\phi \partial\Gamma^\alpha\phi\right\|_{L^1(D_2)}}_{l_{31}}+ \underbrace{\left\|\Omega\Gamma^\beta\phi\partial^2\Gamma^\gamma\phi \Gamma^\alpha\phi\right\|_{L^1(D_2)}}_{l_{32}}\\
&+ \underbrace{\left\|\Omega\Gamma^\beta\phi\partial^2\Gamma^\gamma\phi \langle t-r\rangle(\partial_t-\partial_r)\Gamma^\alpha\phi\right\|_{L^1(D_2)}}_{l_{33}}\\
&+ \underbrace{\left\|\Omega\Gamma^\beta\phi\partial^2\Gamma^\gamma\phi \langle t+r\rangle(\partial_t+\partial_r)\Gamma^\alpha\phi\right\|_{L^1(D_2)}}_{l_{34}}
\end{split}
\end{equation} 
and
\begin{equation}\label{524}
\begin{split}
l_4\lesssim&\underbrace{\left\|\partial\Gamma^\beta\phi\partial\Omega\Gamma^\gamma\phi \partial\Gamma^\alpha\phi\right\|_{L^1(D_2)}}_{l_{41}}+\underbrace{\left\|\partial\Gamma^\beta\phi\partial\Omega\Gamma^\gamma\phi \Gamma^\alpha\phi\right\|_{L^1(D_2)}}_{l_{42}}\\
&+\underbrace{\left\|\partial\Gamma^\beta\phi\partial\Omega\Gamma^\gamma\phi \langle t-r\rangle(\partial_t-\partial_r)\Gamma^\alpha\phi\right\|_{L^1(D_2)}}_{l_{43}}\\
&+\underbrace{\left\|\partial\Gamma^\beta\phi\partial\Omega\Gamma^\gamma\phi \langle t+r\rangle (\partial_t+\partial_r)\Gamma^\alpha\phi\right\|_{L^1(D_2)}}_{l_{44}}.
\end{split}
\end{equation} 
The terms on the right hand sides of \eqref{523} and \eqref{524} will be estimated one by one. First,
\begin{equation}
l_{31}\lesssim\left\{\begin{aligned}
&\left\|\langle t-r\rangle^{-1}\Omega\Gamma^\beta\phi\right\|_{L^\infty}\left\|\langle t-r\rangle\partial^2\Gamma^\gamma\phi\right\|_{L^2}\left\|\partial\Gamma^\alpha\phi\right\|_{L^2}, &|\beta|\leq|\frac{\alpha}{2}|,\\
&\left\|\frac{1}{\langle r\rangle^{1/2}}\Omega\Gamma^\beta\phi\right\|_{L^2}\left\|\langle r\rangle^{1/2}\partial^2\Gamma^\gamma\phi\right\|_{L^\infty}\left\|\partial\Gamma^\alpha\phi\right\|_{L^2}, &|\gamma|\leq|\frac{\alpha}{2}|,
\end{aligned}
\right.
\end{equation}
which implies that
\begin{equation}\label{526}
l_{31}\lesssim E_\mu^{1/2}(\phi)\mathcal{X}_\kappa(\phi)E_\kappa^{1/2}(\phi)+M_\kappa^{1/2}(\phi)E_\mu^{1/2}(\phi)E_\kappa^{1/2}(\phi).
\end{equation}
Since $\langle t\rangle \sim \langle r\rangle$ in $D_2$, then Lemma \ref{lemhardy} implies
\begin{equation}
l_{32}\lesssim\left\{\begin{aligned}
&\left\|\langle r\rangle^{1/2}\langle t-r\rangle^{-1}\Omega\Gamma^\beta\phi\right\|_{L^\infty}\left\|\langle t-r\rangle\partial^2\Gamma^\gamma\phi\right\|_{L^2}\left\|\frac{1}{\langle t\rangle^{1/2}}\Gamma^\alpha\phi\right\|_{L^2}, &|\beta|\leq|\frac{\alpha}{2}|,\\
&\left\|\frac{1}{\langle r\rangle^{1/2}}\Omega\Gamma^\beta\phi\right\|_{L^2}\left\|\langle r\rangle^{1/2}\langle t-r\rangle\partial^2\Gamma^\gamma\phi\right\|_{L^\infty}\left\|\langle t-r\rangle^{-1}\Gamma^\alpha\phi\right\|_{L^2}, &|\gamma|\leq|\frac{\alpha}{2}|,
\end{aligned}
\right.
\end{equation}
which further yields
\begin{equation}\label{528}
l_{32}\lesssim E_\mu^{1/2}(\phi)\mathcal{X}_\kappa(\phi)M_\kappa^{1/2}(\phi)+M_\kappa^{1/2}(\phi)\mathcal{X}_\mu(\phi)E_\kappa^{1/2}(\phi).
\end{equation}
Similarly, we estimate $l_{33}$ as
\begin{equation}
l_{33}\lesssim\left\{\begin{aligned}
&\left\|\langle r\rangle^{1/2}\langle t-r\rangle^{-1}\Omega\Gamma^\beta\phi\right\|_{L^\infty}\left\|\langle t-r\rangle\partial^2\Gamma^\gamma\phi\right\|_{L^2}\left\|\langle t-r\rangle^{1/2}(\partial_t-\partial_r)\Gamma^\alpha\phi\right\|_{L^2}, &|\beta|\leq|\frac{\alpha}{2}|,\\
&\left\|\frac{1}{\langle r\rangle^{1/2}}\Omega\Gamma^\beta\phi\big\|_{L^2}\|\langle r\rangle^{1/2}\langle t-r\rangle^{1/2}\partial^2\Gamma^\gamma\phi\right\|_{L^\infty}\left\|\langle t-r\rangle^{1/2}(\partial_t-\partial_r)\Gamma^\alpha\phi\right\|_{L^2}, &|\gamma|\leq|\frac{\alpha}{2}|,
\end{aligned}
\right.
\end{equation}
which yields
\begin{equation}\label{530}
l_{33}\lesssim E_\mu^{1/2}(\phi)\mathcal{X}_\kappa(\phi)M_\kappa^{1/2}(\phi)+\mathcal{X}_\kappa(\phi)M_\kappa(\phi).
\end{equation}
And for $l_{34}$, noting that $|\beta|+1\leq|\alpha|$, we have
\begin{equation}
l_{34}\lesssim\left\{\begin{aligned}
&\left\|\langle r\rangle^{1/2}\langle t-r\rangle^{-1}\Omega\Gamma^\beta\phi\right\|_{L^\infty}\left\|\langle t-r\rangle\partial^2\Gamma^\gamma\phi\right\|_{L^2}\left\|\langle t+r\rangle^{1/2}(\partial_t+\partial_r)\Gamma^\alpha\phi\right\|_{L^2}, &|\beta|\leq|\frac{\alpha}{2}|,\\
&\left\|\langle t-r\rangle^{-1}\Gamma^{\beta+1}\phi\right\|_{L^2}\left\|\langle r\rangle^{1/2}\langle t-r\rangle\partial^2\Gamma^\gamma\phi\right\|_{L^\infty}\left\|\langle t+r\rangle^{1/2}(\partial_t+\partial_r)\Gamma^\alpha\phi\right\|_{L^2}, &|\gamma|\leq|\frac{\alpha}{2}|,
\end{aligned}
\right.
\end{equation}
which means that
\begin{equation}\label{532}
l_{34}\lesssim E_\mu^{1/2}(\phi)\mathcal{X}_\kappa(\phi)M_\kappa^{1/2}(\phi)+E_\kappa^{1/2}(\phi)\mathcal{X}_\mu(\phi)M_\kappa^{1/2}(\phi).
\end{equation}
Inserting \eqref{526}, \eqref{528}, \eqref{530} and \eqref{532} into \eqref{523} and combining Lemma \ref{lem3.3}, we get
\begin{equation}\label{533}
l_3\lesssim E_\mu^{1/2}(\phi)M_\kappa(\phi).
\end{equation}
Since 
\begin{equation}
l_{41}\lesssim\left\{\begin{aligned}
&\left\|\partial\Gamma^\beta\phi\right\|_{L^\infty}\left\|\partial\Omega\Gamma^\gamma\phi\right\|_{L^2}\left\|\partial\Gamma^\alpha\phi\right\|_{L^2}, &|\beta|\leq|\frac{\alpha}{2}|,\\
&\left\|\partial\Gamma^\beta\phi\right\|_{L^2}\left\|\partial\Omega\Gamma^\gamma\phi\right\|_{L^\infty}\left\|\partial\Gamma^\alpha\phi\right\|_{L^2}, &|\gamma|\leq|\frac{\alpha}{2}|,
\end{aligned}
\right.
\end{equation}

\begin{equation}
l_{42}\lesssim\left\{\begin{aligned}
&\left\|\langle r\rangle^{1/2}\partial\Gamma^\beta\phi\right\|_{L^\infty}\left\|\partial\Omega\Gamma^\gamma\phi\right\|_{L^2}\left\|\frac{1}{\langle r\rangle^{1/2}}\Gamma^\alpha\phi\right\|_{L^2}, &|\beta|\leq|\frac{\alpha}{2}|,\\
&\left\|\partial\Gamma^\beta\phi\right\|_{L^2}\left\|\langle r\rangle^{1/2}\partial\Omega\Gamma^\gamma\phi\right\|_{L^\infty}\left\|\frac{1}{\langle r\rangle^{1/2}}\Gamma^\alpha\phi\right\|_{L^2}, &|\gamma|\leq|\frac{\alpha}{2}|,
\end{aligned}
\right.
\end{equation}

\begin{equation}
l_{43}\lesssim\left\{\begin{aligned}
&\left\|\langle r\rangle^{1/2}\partial\Gamma^\beta\phi\right\|_{L^\infty}\left\|\partial\Omega\Gamma^\gamma\phi\right\|_{L^2}\left\|\langle t-r\rangle^{1/2}(\partial_t-\partial_r)\Gamma^\alpha\phi\right\|_{L^2}, &|\beta|\leq|\frac{\alpha}{2}|,\\
&\left\|\partial\Gamma^\beta\phi\right\|_{L^2}\left\|\langle r\rangle^{1/2}\partial\Omega\Gamma^\gamma\phi\right\|_{L^\infty}\left\|\langle t-r\rangle^{1/2}(\partial_t-\partial_r)\Gamma^\alpha\phi\right\|_{L^2}, &|\gamma|\leq|\frac{\alpha}{2}|,
\end{aligned}
\right.
\end{equation}
and 
\begin{equation}
l_{44}\lesssim\left\{\begin{aligned}
&\left\|\langle r\rangle^{1/2}\partial\Gamma^\beta\phi\right\|_{L^\infty}\left\|\partial\Omega\Gamma^\gamma\phi\right\|_{L^2}\left\|\langle t+r\rangle^{1/2}(\partial_t+\partial_r)\Gamma^\alpha\phi\right\|_{L^2}, &|\beta|\leq|\frac{\alpha}{2}|,\\
&\left\|\partial\Gamma^\beta\phi\right\|_{L^2}\left\|\langle r\rangle^{1/2}\partial\Omega\Gamma^\gamma\phi\right\|_{L^\infty}\left\|\langle t+r\rangle^{1/2}(\partial_t+\partial_r)\Gamma^\alpha\phi\right\|_{L^2}, &|\gamma|\leq|\frac{\alpha}{2}|,
\end{aligned}
\right.
\end{equation}
we also have
\begin{equation}\label{538}
l_4\leq l_{41}+l_{42}+l_{43}+l_{44}\lesssim E_\mu^{1/2}(\phi)M_\kappa(\phi).
\end{equation}
Finally, combing \eqref{520}, \eqref{521}-\eqref{522}, \eqref{533} and \eqref{538}, we obtain
\begin{equation}\label{539}
L_{11}(D_2)+L_{2}(D_2)+ L_{3}(D_2) \lesssim \langle t\rangle^{-1} E_\mu^{1/2}(\phi)M_\kappa(\phi).
\end{equation}
As the term $L_{12}$ in \eqref{4.5} is a certain kind of $L_2$ and $L_3$, combining \eqref{518} and \eqref{539}, we get the following result
\begin{equation}
  \begin{split}
        \frac{d}{dt}\bar{M}_\kappa( \phi(t))\lesssim\langle t\rangle^{-1}E_{\mu}^{1/2}( \phi(t))M_\kappa( \phi(t)), \qquad 0<t<T.
  \end{split}
\end{equation}

\subsection{General energy estimate}
Let $|\alpha|\leq \mu-1$ and $\mu=\kappa-1$. Multiplying $\partial_t\Gamma^\alpha u^i$ on both sides of equations \eqref{Gphieq}, integrating on $\mathbb{R}^2$, we have
\begin{equation}
  \begin{split}
    \frac{d}{dt}E_{\mu}(\phi)=&-\frac12\sum_{|\alpha|=0}^{\mu-1}\frac{d}{dt}\int_{\mathbb{R}^2} \left(2g^{ijk}\partial_k\phi +3h^{ijkl}\partial _k\phi\partial_l\phi\right)\partial_{j}\Gamma^\alpha\phi \partial_i\Gamma^\alpha\phi dx\\
&+\frac12\sum_{|\alpha|=0}^{\mu-1}\int_{\mathbb{R}^2} \partial_t\left(2g^{ijk}\partial_k\phi +3h^{ijkl}\partial _k\phi\partial_l\phi\right)\partial_{j}\Gamma^\alpha\phi \partial_i\Gamma^\alpha\phi dx\\
    &+\sum_{|\alpha|=0}^{\mu-1}\sum_{|\beta+\gamma|\leq|\alpha|\atop |\beta|, |\gamma|<|\alpha|}C_\alpha^\beta\int_{\mathbb{R}^2} g^{ijk}\partial_i\left(\partial_j\Gamma^\beta\phi \partial_k\Gamma^\gamma\phi\right) \partial_t\Gamma^\alpha\phi dx\\
    &+\sum_{|\alpha|=0}^{\mu-1}\sum_{|\beta+\gamma+\delta|\leq|\alpha|\atop |\beta|, |\gamma|, |\delta|<|\alpha|}C_\alpha^\beta C_{\gamma+\delta}^\gamma\int_{\mathbb{R}^2} h^{ijkl}\partial _i\left(\partial_j\Gamma^\beta\phi\partial_k\Gamma^\gamma\phi\partial_l\Gamma^\delta\phi\right) \partial_t\Gamma^\alpha\phi dx.
  \end{split}
\end{equation}
Let
\begin{equation}
  \tilde{E}_\mu(\phi)=E_\mu(\phi)+\frac12\sum_{|\alpha|=0}^{\mu-1}\int_{\mathbb{R}^2} \left(2g^{ijk}\partial_k\phi +3h^{ijkl}\partial _k\phi\partial_l\phi\right)\partial_{j}\Gamma^\alpha\phi \partial_i\Gamma^\alpha\phi dx.
  \end{equation}
Due to the fact that
\begin{equation}
\int_{\mathbb{R}^2}\left (2g^{ijk}\partial_k\phi +3h^{ijkl}\partial _k\phi\partial_l\phi\right)\partial_{j}\Gamma^\alpha\phi \partial_i\Gamma^\alpha\phi dx
\lesssim \left\|\partial \phi+(\partial\phi)^2\right\|_{L^\infty}\left\|\partial\Gamma^\alpha\phi\right\|_{L^2}^2
\lesssim \epsilon E_\mu(\phi)
\end{equation}
for $\epsilon\ll 1$, it holds that
\begin{equation}
\tilde{E}_\mu(\phi)\approx E_\mu(\phi).
\end{equation}
So we have
\begin{equation}
  \begin{split}
    \frac{d}{dt}\tilde{E}_{\mu}(\phi)=
&\frac12\sum_{|\alpha|=0}^{\mu-1}\underbrace{\int_{\mathbb{R}^2} \partial_t\left(2g^{ijk}\partial_k\phi +3h^{ijkl}\partial _k\phi\partial_l\phi\right)\partial_{j}\Gamma^\alpha\phi \partial_i\Gamma^\alpha\phi dx}_{G_1}\\
    &+\sum_{|\alpha|=0}^{\mu-1}\sum_{|\beta+\gamma|\leq|\alpha|\atop |\beta|, |\gamma|<|\alpha|}C_\alpha^\beta\underbrace{\int_{\mathbb{R}^2} g^{ijk}\partial_i\left(\partial_j\Gamma^\beta\phi \partial_k\Gamma^\gamma\phi\right) \partial_t\Gamma^\alpha\phi dx}_{G_2}\\
    &+\sum_{|\alpha|=0}^{\mu-1}\sum_{|\beta+\gamma+\delta|\leq|\alpha|\atop |\beta|, |\gamma|, |\delta|<|\alpha|}C_\alpha^\beta C_{\gamma+\delta}^\gamma\underbrace{\int_{\mathbb{R}^2} h^{ijkl}\partial _i\left(\partial_j\Gamma^\beta\phi\partial_k\Gamma^\gamma\phi\partial_l\Gamma^\delta\phi\right) \partial_t\Gamma^\alpha\phi dx}_{G_3}.
  \end{split}
\end{equation}

For $G_1$, the null conditions help us gain more time decay. Using Lemma \ref{nulllem}, we then get
\begin{equation}\label{546}
  \begin{split}
    G_1\lesssim
\langle t\rangle^{-1}\Bigg[&\underbrace{\left\|\left(\partial\Gamma\phi +\partial \phi\partial\Gamma\phi+\Gamma\phi\partial^2\phi\right)\left(\partial\Gamma^\alpha\phi\right)^2\right\|_{L^1}}_{G_{11}}+\underbrace{\left\|\left(\partial^2\phi +\partial \phi\partial^2\phi\right)\partial\Gamma^\alpha\phi\Gamma^{\alpha+1}\phi\right\|_{L^1}}_{G_{12}}\\
&+\underbrace{\left\|\langle t-r\rangle\left(\partial^2\phi +\partial \phi\partial^2\phi\right)\left(\partial\Gamma^\alpha\phi\right)^2\right\|_{L^1}}_{G_{13}}\Bigg].
\end{split}
\end{equation}
Noting that 
\begin{equation}
1\lesssim \frac{\langle r\rangle^{1/2}\langle t-r\rangle^{1/2}}{\langle t\rangle^{1/2}},
\end{equation}
then the first term on the right hand side of \eqref{546} is controlled by
\begin{equation}\label{548}
\begin{split}
G_{11}\lesssim \langle t\rangle^{-1/2}\Bigg[ &\left\|\langle r\rangle^{1/2}\langle t-r\rangle^{1/2}\partial\Gamma\phi (\partial\Gamma^\alpha\phi)^2\right\|_{L^1}+\left\|\langle r\rangle^{1/2}\langle t-r\rangle^{1/2}\partial\phi\partial\Gamma\phi (\partial\Gamma^\alpha\phi)^2\right\|_{L^1}\\
&+ \left\|\langle r\rangle^{1/2}\langle t-r\rangle^{1/2}\Gamma\phi \partial^2\phi(\partial\Gamma^\alpha\phi)^2\right\|_{L^1}\Bigg]\\
\qquad\lesssim \langle t\rangle^{-1/2}\Bigg[ &\left\|\langle r\rangle^{1/2}\langle t-r\rangle^{1/2}\partial\Gamma\phi \right\|_{L^\infty}\left\|\partial\Gamma^\alpha\phi\right\|_{L^2}^2\\
&+\left\|\langle r\rangle^{1/2}\langle t-r\rangle^{1/2}\partial\phi\right\|_{L^\infty}\left\|\partial\Gamma\phi \right\|_{L^\infty}\left\|\partial\Gamma^\alpha\phi\right\|_{L^2}^2\\
&+\left \|\langle r\rangle^{1/2}\langle t-r\rangle^{-1}\Gamma\phi\right\|_{L^\infty}\left\|\langle t-r\rangle \partial^2\phi\right\|_{L^\infty}\left\|\langle t-r\rangle^{1/2}\partial\Gamma^\alpha\phi\right\|_{L^2}\left\|\partial\Gamma^\alpha\phi\right\|_{L^2}\Bigg]\\
\qquad\lesssim\langle t\rangle^{-1/2}\Big[ &\mathcal{X}_\mu(\phi)E_\mu(\phi)+\varepsilon\mathcal{X}_\mu(\phi)E_\mu(\phi)+\varepsilon\mathcal{X}_\mu(\phi)M_\mu^{1/2}(\phi)E_\mu^{1/2}(\phi)\Big]\\
 \lesssim\langle t\rangle^{-1/2}&E_\mu(\phi)M_\mu^{1/2}(\phi),
\end{split}
\end{equation}
where we have used Lemma \ref{lem0}, Lemma \ref{lemhardy}, Lemma \ref{lemM1} and Lemma \ref{lem3.3}. For $G_{12}$, we split the whole space into two parts: $D_1=\{x: |x|\leq\frac{t}{2}\ \text {or}\ 2t\leq|x|\leq t+1\}$ and $D_2=\{x:\frac{t}2<|x|<2t\}$.  In domain $D_1$, noting that
\begin{equation}
1\lesssim\frac{\langle t-r\rangle^{1/2}}{\langle t\rangle^{1/2}},
\end{equation}
it holds
\begin{equation}
\begin{split}
G_{12}(D_1)\lesssim \langle t\rangle^{-1/2}\Bigg[ &\left\|\langle t-r\rangle^{1/2}\partial^2\phi \partial\Gamma^\alpha\phi\Gamma^{\alpha+1}\phi\right\|_{L^1}+\left\|\langle t-r\rangle^{1/2}\partial\phi\partial^2\phi \partial\Gamma^\alpha\phi\Gamma^{\alpha+1}\phi\right\|_{L^1}\Bigg]\\
\qquad\lesssim \langle t\rangle^{-1/2}\Bigg[ &\left\|\langle t-r\rangle\partial^2\phi \right\|_{L^\infty}\left\|\partial\Gamma^\alpha\phi\right\|_{L^2}\left\|\frac{1}{\langle t\rangle^{1/2}}\Gamma^{\alpha+1}\phi\right\|_{L^2}\\
&+\|\partial\phi\|_{L^\infty}\left\|\langle t-r\rangle\partial^2\phi \right\|_{L^\infty}\left\|\partial\Gamma^\alpha\phi\right\|_{L^2}\left\|\frac{1}{\langle t\rangle^{1/2}}\Gamma^{\alpha+1}\phi\right\|_{L^2}\Bigg]\\
\qquad\lesssim\langle t\rangle^{-1/2}\Big[ &\mathcal{X}_\mu(\phi)E_\mu^{1/2}(\phi)M_\kappa^{1/2}(\phi)+\varepsilon\mathcal{X}_\mu(\phi)E_\mu^{1/2}(\phi)M_\kappa^{1/2}(\phi)\Big]\\
\lesssim\langle t\rangle^{-1/2}&E_\mu(\phi)M_\kappa^{1/2}(\phi).
\end{split}
\end{equation}
And in domain $D_2$, since $t\sim r$, we then get
\begin{equation}
\begin{split}
G_{12}(D_2)\lesssim \langle t\rangle^{-1/2}\Bigg[ &\left\|\langle r\rangle^{1/2}\partial^2\phi \partial\Gamma^\alpha\phi\Gamma^{\alpha+1}\phi\right\|_{L^1}+\left\|\langle r\rangle^{1/2}\partial\phi\partial^2\phi \partial\Gamma^\alpha\phi\Gamma^{\alpha+1}\phi\right\|_{L^1}\Bigg]\\
\qquad\lesssim \langle t\rangle^{-1/2}\Bigg[ &\left\|\langle r\rangle^{1/2}\langle t-r\rangle\partial^2\phi \right\|_{L^\infty}\left\|\partial\Gamma^\alpha\phi\right\|_{L^2}\left\|\frac{1}{\langle t-r\rangle}\Gamma^{\alpha+1}\phi\right\|_{L^2}\\
&+\left\|\partial\phi\right\|_{L^\infty}\left\|\langle r\rangle^{1/2}\langle t-r\rangle\partial^2\phi \right\|_{L^\infty}\left\|\partial\Gamma^\alpha\phi\right\|_{L^2}\left\|\frac{1}{\langle t-r\rangle}\Gamma^{\alpha+1}\phi\right\|_{L^2}\Bigg]\\
\qquad\lesssim\langle t\rangle^{-1/2}\Big[ &\mathcal{X}_\mu(\phi)E_\mu^{1/2}(\phi)E_\kappa^{1/2}(\phi)+\varepsilon\mathcal{X}_\mu(\phi)E_\mu^{1/2}(\phi)E_\kappa^{1/2}(\phi)\Big]\\
\lesssim\langle t\rangle^{-1/2}&E_\mu(\phi)E_\kappa^{1/2}(\phi).
\end{split}
\end{equation}
Combing these two cases in $D_1$ and $D_2$ together, we arrive at
\begin{equation}\label{552}
G_{12}\lesssim\langle t\rangle^{-1/2}E_\mu(\phi)M_\kappa^{1/2}(\phi).
\end{equation}
The estimate for the last term ($G_{13}$) on the right hand side of \eqref{546} is given as follows
\begin{equation}\label{553}
\begin{split}
G_{13}\lesssim \langle t\rangle^{-1/2}\Bigg[ &\left\|\langle r\rangle^{1/2}\langle t-r\rangle^{3/2}\partial^2\phi (\partial\Gamma^\alpha\phi)^2\right\|_{L^1}+\left\|\langle r\rangle^{1/2}\langle t-r\rangle^{3/2}\partial\phi\partial^2\phi (\partial\Gamma^\alpha\phi)^2\right\|_{L^1}\Bigg]\\
\qquad\lesssim \langle t\rangle^{-1/2}\Bigg[ &\left\|\langle r\rangle^{1/2}\langle t-r\rangle\partial^2\phi \right\|_{L^\infty}\left\|\langle t-r\rangle^{1/2}\partial\Gamma^\alpha\phi\right\|_{L^2}\left\|\partial\Gamma^\alpha\phi\right\|_{L^2}\\
&+\left\|\partial\phi\right\|_{L^\infty}\left\|\langle r\rangle^{1/2}\langle t-r\rangle\partial^2\phi \right\|_{L^\infty}\left\|\langle t-r\rangle^{1/2}\partial\Gamma^\alpha\phi\right\|_{L^2}\|\partial\Gamma^\alpha\phi\|_{L^2}\Bigg]\\
\qquad\lesssim\langle t\rangle^{-1/2}\Big[ &\mathcal{X}_\mu(\phi)M_\mu^{1/2}(\phi)E_\mu^{1/2}(\phi)+\varepsilon\mathcal{X}_\mu(\phi)M_\mu^{1/2}(\phi)E_\mu^{1/2}(\phi)\Big]\\
\ \lesssim\langle t\rangle^{-1/2}&E_\mu(\phi)M_\mu^{1/2}(\phi).
\end{split}
\end{equation}
Hence, inserting \eqref{548}, \eqref{552} and \eqref{553} into \eqref{546}, we obtain
\begin{equation}\label{554}
G_1\lesssim\langle t\rangle^{-3/2}E_\mu(\phi)M_\kappa^{1/2}(\phi).
\end{equation}

For $G_2$, due to the null conditions, it holds
\begin{equation}\label{555}
\begin{split}
    G_2\lesssim
\langle t\rangle^{-1}\Bigg[&\underbrace{\left\|\partial\Gamma^{\beta+1}\phi \partial\Gamma^\gamma\phi\partial\Gamma^\alpha\phi\right\|_{L^1}}_{G_{21}}+\underbrace{\left\|\partial^2\Gamma^{\beta}\phi \Gamma^{\gamma+1}\phi\partial\Gamma^\alpha\phi\right\|_{L^1}}_{G_{22}}\\
&+\underbrace{\left\|\langle t-r\rangle\partial^2\Gamma^{\beta}\phi \partial\Gamma^\gamma\phi\partial\Gamma^\alpha\phi\right\|_{L^1}}_{G_{23}}\Bigg].
  \end{split}
\end{equation}
The terms on the right hand side of \eqref{555} can be controlled as
\begin{equation}
\begin{aligned}
G_{21}&\lesssim\left\{\begin{aligned}
\langle t\rangle^{-1/2}\left\|\langle r\rangle^{1/2}\langle t-r\rangle^{1/2}\partial\Gamma^{\beta+1}\phi\right\|_{L^\infty}\left\| \partial\Gamma^\gamma\phi\right\|_{L^2}\left\|\partial\Gamma^\alpha\phi\right\|_{L^2},&\ |\beta|\leq \left|\frac{\alpha}{2}\right|\\
\langle t\rangle^{-1/2}\left\|\partial\Gamma^{\beta+1}\phi\right\|_{L^2}\left\|\langle r\rangle^{1/2}\langle t-r\rangle^{1/2} \partial\Gamma^\gamma\phi\right\|_{L^\infty}\left\|\partial\Gamma^\alpha\phi\right\|_{L^2},&\ |\gamma|\leq \left|\frac{\alpha}{2}\right|
\end{aligned}\right.\\
&\lesssim\langle t\rangle^{-1/2}\mathcal{X}_\mu(\phi)E_\mu(\phi)\\
&\lesssim\langle t\rangle^{-1/2}E_\mu^{3/2}(\phi),
\end{aligned}
\end{equation}

\begin{equation}
\begin{aligned}
G_{22}&\lesssim\left\{\begin{aligned}
\langle t\rangle^{-1/2}\left\|\langle r\rangle^{1/2}\langle t-r\rangle\partial^2\Gamma^{\beta}\phi\right\|_{L^\infty}\left\|\langle t-r\rangle^{-1} \Gamma^{\gamma+1}\phi\right\|_{L^2}\left\|\langle t-r\rangle^{1/2}\partial\Gamma^\alpha\phi\right\|_{L^2},&\ |\beta|\leq \left|\frac{\alpha}{2}\right|\\
\langle t\rangle^{-1/2}\left\|\langle t-r\rangle\partial^2\Gamma^{\beta}\phi\right\|_{L^2}\left\|\langle r\rangle^{1/2}\langle t-r\rangle^{-1} \Gamma^{\gamma+1}\phi\right\|_{L^\infty}\left\|\langle t-r\rangle^{1/2}\partial\Gamma^\alpha\phi\right\|_{L^2},&\ |\gamma|\leq \left|\frac{\alpha}{2}\right|
\end{aligned}\right.\\
&\lesssim\langle t\rangle^{-1/2}\mathcal{X}_\mu(\phi)E_\mu^{1/2}(\phi)M_\mu^{1/2}(\phi)\\
&\lesssim\langle t\rangle^{-1/2}E_\mu(\phi)M_\mu^{1/2}(\phi),
\end{aligned}
\end{equation}
and 
\begin{equation}\label{558}
\begin{aligned}
G_{23}&\lesssim\left\{\begin{aligned}
\langle t\rangle^{-1/2}\left\|\langle r\rangle^{1/2}\langle t-r\rangle\partial^2\Gamma^{\beta}\phi\right\|_{L^\infty}\left\|\langle t-r\rangle^{1/2} \partial\Gamma^\gamma\phi\right\|_{L^2}\left\|\partial\Gamma^\alpha\phi\right\|_{L^2},&\ |\beta|\leq \left|\frac{\alpha}{2}\right|\\
\langle t\rangle^{-1/2}\left\|\langle t-r\rangle\partial^2\Gamma^{\beta}\phi\right\|_{L^2}\left\|\langle r\rangle^{1/2}\langle t-r\rangle^{1/2} \partial\Gamma^\gamma\phi\right\|_{L^\infty}\left\|\partial\Gamma^\alpha\phi\right\|_{L^2},&\ |\gamma|\leq \left|\frac{\alpha}{2}\right|
\end{aligned}\right.\\
&\lesssim\langle t\rangle^{-1/2}\mathcal{X}_\mu(\phi)E_\mu^{1/2}(\phi)M_\mu^{1/2}(\phi)\\
&\lesssim\langle t\rangle^{-1/2}E_\mu(\phi)M_\mu^{1/2}(\phi).
\end{aligned}
\end{equation}
Putting \eqref{555}-\eqref{558} together, we get
\begin{equation}\label{559}
G_2\lesssim\langle t\rangle^{-3/2}E_\mu(\phi)M_\mu^{1/2}(\phi).
\end{equation}

The term $G_3$ can be estimated similarly as $G_2$, by taking another $L^\infty$ norm for the function with lower-order derivatives.  We omit the details and write out the conclusion 
\begin{equation}\label{560}
G_3\lesssim\langle t\rangle^{-3/2}\varepsilon E_\mu(\phi)M_\mu^{1/2}(\phi),
\end{equation}
and finally we come to
\begin{equation}
 \frac{d}{dt}\tilde{E}_{\mu}(\phi)\lesssim\langle t\rangle^{-3/2}E_\mu(\phi)M_\kappa^{1/2}(\phi).
\end{equation}

\section{Application to elastic waves}\label{sec6}
Let  $\Psi:\mathbb{R}\times\mathbb{R}^2\to\mathbb{R}^2$ with $\Psi=\Psi(t,x)$ be the motion of an elasticity, satisfying $\Psi(0,x)=x=(x_1,x_2)$.
For homogeneous, isotropic hyperelastic materials with constant unit density, the equations of motion can be derived through application of the principle of least action to $\mathcal{L}$:
\begin{align*}
\mathcal{L}:=\int\int_{\mathbb{R}^2}\Big(\frac12|\partial_t\Psi|^2-W(FF^T)\Big)dxdt,
\end{align*}
where  $W$ is the stored energy depending on the deformation gradient $F_{ij}=\partial_{x_j}\Psi_i$. Then the  equations take the form:
\begin{equation}\label{el}
  \frac{\partial^2\Psi_i}{\partial t^2 }-\frac{\partial}{\partial{x_l}}\frac{\partial W(FF^T)}{\partial F_{il}}=0.
\end{equation}
Let $U:=\Psi-x$ be the displacement.
Using the Taylor expansion of $W$ near equilibrium state, the motion of displacement $U$ satisfies a quasilinear wave system with multiple wave-speeds:
\begin{equation} \label{waves}
    \partial_t^2 U-c_2^2\Delta U-(c_1^2-c_2^2)\nabla(\nabla\cdot U)=N(\nabla U,\nabla^2U),
\end{equation}
where $N(\nabla U,\nabla^2U)$ is composed of $\nabla U\nabla^2U$-form, constants $c_1>c_2>0$ are presented as the pressure wave speed and shear wave speed respectively.

\subsection{Admissible harmonic elastodynamics}
As discussed by John \cite{John60,John66}, an elastic material is harmonic, if the pseudo-irrotationality of its deformation at the initial time can be always preserved at any time. That is,
\begin{equation}
    F_{ij}=F_{ji}\quad\text{or}\quad \partial_{x_j}\Psi_i=\partial_{x_i}\Psi_j,\quad \forall \quad t \geq0.
\end{equation}
Thus, system \eqref{waves} reduces to a scalar wave equation. 

In this paper, we consider the admissible harmonic elastic materials proposed by An-Chen-Yu \cite{ACY}:
\begin{equation} \label{hw}
    \partial_t^2 U-c_2^2\Delta U-(c_1^2-c_2^2)\nabla(\nabla\cdot U)=\text{grad}\ G(\partial U),
\end{equation}
where $G(\partial U)$ is a scalar function of $\partial U$.
The initial data of \eqref{hw} is given by:
\begin{equation}\label{hwdata}
    U(0,x)= U_0(x),\ U_t(0,x)= U_1(x),
\end{equation}
where $U_0(x),U_1(x)\in C^\infty(\mathbb{R}^2)$ with compact support $\{x: |x|\leq 1^-\}$.

Expanding $G$ into its Taylor series around $\partial U=0$, we obtain
\begin{equation}\label{68}
    G(\partial U)=G(0)+s_{jk}\partial_jU_k+g_{jk,mn}\partial_jU_k\partial_mU_n+h_{jk,mn}^{ab}\partial_jU_k\partial_mU_n\partial_aU_b+\cdots.
\end{equation}
Here repeated indices indicate summation from 1 to 2. 
Without loss of generality, we suppose that the admissible harmonic elasticity satisfies 
\begin{equation}
    s_{jk}=0.
\end{equation}
Hence equation \eqref{hw} is a quasilinear wave system in the form of 
\begin{equation} \label{6.7}
    LU=\text{grad}(g_{jk,mn}\partial_jU_k\partial_mU_n+h_{jk,mn}^{ab}\partial_jU_k\partial_mU_n\partial_aU_b+\cdots),
\end{equation}
where 
$$LU=\partial_t^2 U-c_2^2\Delta U-(c_1^2-c_2^2)\nabla(\nabla\cdot U).$$
With some kinds of null conditions and small initial data, we only keep the quadratic and cubic nonlinear terms in \eqref{6.7}.

\subsection{Null conditions}
The main purpose of this section is to show the global well-posedness of the Cauchy problem of certain elastic waves with small initial data. As we have known, the null conditions for quadratic and cubic quasilinear terms are necessary. 

As in Section 2.1 of Li \cite{Li}, system \eqref{6.7} satisfies the {\it null condition} if each small plane wave solution $U=U(ct+x\cdot w)$ $(U(0)=0, U'(0)=0)$ to the linearized system $LU=0$ is always a solution to the nonlinear system \eqref{6.7}. Here $c$ is a constant. For $LU=0$, there are two families of planer waves:
\begin{equation}
    \begin{split}
    U^{rad}=& w\psi_1(c_1t+r),\ \forall\ \psi_1\ \text{is a scalar function}, \\
    U^{tan}=& w^\bot\psi_2(c_2t+r), \ \forall\ \psi_2\ \text{is a scalar function},
    \end{split}
\end{equation}
see Sideris \cite{S2}.
We say \eqref{6.7} satisfies {\it the first null condition} if $$\nabla\left[g_{jk,mn}\partial_jU_k\partial_mU_n\right]=0,\ \forall\ U\in \{U^{rad}, U^{tan}\}$$
and {\it the second null condition} if
$$\nabla\left[h_{jk,mn}^{ab}\partial_jU_k\partial_mU_n\partial_aU_b\right]=0, \ \forall\ U\in \{U^{rad}, U^{tan}\}.$$

We consider a kind of the so-called admissible harmonic elasticity, whose dynamics admits:
\begin{equation} \label{hweq}
\begin{split}
    &\partial_t^2 U-c_2^2\Delta U-(c_1^2-c_2^2)\nabla(\nabla\cdot U)\\
    =&\nabla\left[g_0(\nabla\cdot U)^2+g_1(\nabla\cdot U)(\nabla^\bot\cdot U)+g_2(\nabla^\bot\cdot U)^2\right]\\
    &+\nabla\left[h_0(\nabla\cdot U)^3+h_1(\nabla\cdot U)^2(\nabla^\bot\cdot U)+h_2(\nabla\cdot U)(\nabla^\bot\cdot U)^2+h_3(\nabla^\bot\cdot U)^3\right].
    \end{split}
\end{equation}
It is easy to verify that:
\begin{lem}
  If $g_0=0$ and $h_0=0$, then system \eqref{hweq} satisfies the first and second null conditions, respectively.
\end{lem}

For system \eqref{hweq}, we have the following theorem, a detailed version of which will be given in Theorem \ref{thm62} at the end of this section.
\begin{thm}\label{th61}
 If $g_0=h_0=0$, the Cauchy problem of system \eqref{hweq} with small initial data \eqref{hwdata} admits a unique almost global classical solution. If $g_0=g_2=h_0=0$, the Cauchy problem of system \eqref{hweq} with small initial data \eqref{hwdata} admits a unique global classical solution. 
\end{thm}
\begin{remark}
 For the global existence result in Theorem \ref{th61}, in addition to the null conditions, an additional algebraic condition $g_2=0$ is also posed. Nevertheless, it is still a quasilinear wave system with multiple speeds and quadratic nonlinear terms.
\end{remark}

\subsection{Notations and additional prelimilaries}
Using Helmholtz decomposition, the vector field $U$ can be projected into two parts, i.e., the divergence-free part $U^{df}$ and curl-free part $U^{cf}$: 
\begin{equation}\label{611}
    U=-\nabla(-\Delta)^{-1}\nabla\cdot U+\nabla^\bot(-\Delta)^{-1}\nabla^\bot\cdot U=U^{cf}+U^{df}
\end{equation}
satisfying 
\begin{equation}
    \nabla\cdot U^{df}=0,\qquad \nabla^\bot\cdot U^{cf}=0.
\end{equation}
If the null conditions hold, i.e. $g_0=h_0=0$, then system \eqref{hweq} can be reduced to
\begin{equation}\label{6.13n}
\left\{
    \begin{aligned}
    \partial_t^2 U^{cf}-c_1^2\Delta U^{cf}=&\nabla\left[g_1\left(\nabla\cdot U^{cf}\right)\left(\nabla^\bot\cdot U^{df}\right)+g_2\left(\nabla^\bot\cdot U^{df}\right)^2\right]\\
    &+\nabla\Bigg[h_1\left(\nabla\cdot U^{cf}\right)^2\left(\nabla^\bot\cdot U^{df}\right)+h_2\left(\nabla\cdot U^{cf}\right)\left(\nabla^\bot\cdot U^{df}\right)^2\\
    &+h_3\left(\nabla^\bot\cdot U^{df}\right)^3\Bigg],\\
   \partial_t^2 U^{df}-c_2^2\Delta U^{df}=&0.
\end{aligned}\right.
\end{equation}

For the equation of $U$ in \eqref{hw}, the vector fields $\Gamma$ are modified by $\Gamma=(\partial_t,\partial_1,\partial_2, \tilde{S},\tilde{\Omega})$, where
$$\tilde{S}=S-1,$$
$$\tilde{\Omega}=\Omega I+O, \ O=\left(\begin{matrix}
    0&-1\\
    1&0
\end{matrix}\right),$$
and for the initial data, $\Gamma$ is replaced by $\Gamma_0=\{\partial_1,\partial_2, r\partial_r,\Omega\}$.
Using these vector fields, we have
\begin{equation}\label{n620}
\begin{aligned}
     &\partial_t^2\Gamma^\alpha U-c_2^2\Delta\Gamma^\alpha U-(c_1^2-c_2^2)\nabla(\nabla\cdot\Gamma^\alpha U)\\
     =&\sum_{|\beta+\gamma|\leq|\alpha|}C_\alpha^\beta  
     \nabla\left[g_1\left(\nabla\cdot\Gamma^\beta U\right)\left(\nabla^\bot\cdot \Gamma^\gamma U\right)+g_2\left(\nabla^\bot\cdot\Gamma^\beta U\right)\left(\nabla^\bot\cdot\Gamma^\gamma U\right)\right]\\
    &+\sum_{|\alpha+\beta+\delta|\leq|\alpha|}C_\alpha^\beta C_{\gamma+\delta}^\gamma\nabla\Bigg[h_1\left(\nabla\cdot\Gamma^\beta U\right)\left(\nabla\cdot\Gamma^\gamma U\right)\left(\nabla^\bot\cdot\Gamma^\delta U\right)\\
    &\qquad\qquad\qquad\qquad+h_2\left(\nabla\cdot \Gamma^\beta U\right)\left(\nabla^\bot\cdot \Gamma^\gamma U\right)\left(\nabla^\bot\cdot \Gamma^\delta U\right)\\
    &\qquad\qquad\qquad\qquad+h_3\left(\nabla^\bot\cdot \Gamma^\beta U\right)\left(\nabla^\bot\cdot \Gamma^\gamma U\right)\left(\nabla^\bot\cdot \Gamma^\delta U\right)\Bigg].
\end{aligned}
\end{equation}

The basic general energy associated with the linear part of \eqref{n620} is 
\begin{equation}
    \begin{aligned}
        E_{\kappa}(U)=\frac12\sum_{|\alpha|\leq\kappa-1}\int_{\mathbb{R}^2}\left[|\partial_t\Gamma^\alpha U|^2+c_2^2|\nabla \Gamma^\alpha U|^2+(c_1^2-c_2^2)|\nabla \cdot\Gamma^\alpha U|^2\right]dx.
    \end{aligned}
\end{equation}
\begin{remark}
    It is obvious that 
    \begin{equation}
    \begin{split}
        E_\kappa(U)=&\frac12\sum_{|\alpha|\leq\kappa-1}\int_{\mathbb{R}^2}\left[|\partial_t\Gamma^\alpha U|^2+c_1^2|\nabla \cdot\Gamma^\alpha U|^2+c_2^2|\nabla^\bot\cdot \Gamma^\alpha U|^2\right]dx\\
        =&\underbrace{\frac12\sum_{|\alpha|\leq\kappa-1}\int_{\mathbb{R}^2}\left[|\partial_t\Gamma^\alpha U^{cf}|^2+c_1^2|\nabla \cdot\Gamma^\alpha U^{cf}|^2\right]dx}_{E_\kappa(U^{cf})}\\
        &+\underbrace{\frac12\sum_{|\alpha|\leq\kappa-1}\int_{\mathbb{R}^2}\left[|\partial_t\Gamma^\alpha U^{df}|^2+c_2^2|\nabla^\bot\cdot \Gamma^\alpha U^{df}|^2\right]dx}_{E_\kappa(U^{df})}\\
        =:&E_{\kappa,\kappa}(U^{cf},U^{df}).
        \end{split}
    \end{equation}
\end{remark}
We will also use the weighted $L^2$ energy introduced by Sideris \cite{S1}:
\begin{equation}
    \chi_{\kappa}(U)=\underbrace{\sum_{|\alpha|\leq\kappa-2}\left\|\langle c_1 t-r\rangle\partial\nabla\Gamma^\alpha U^{cf}\right\|_{L^2(\mathbb{R}^2)}}_{\chi_{\kappa}(U^{cf})}+\underbrace{\sum_{|\alpha|\leq\kappa-2}\left\|\langle c_2t-r\rangle\partial\nabla\Gamma^\alpha U^{df}\right\|_{L^2(\mathbb{R}^2)}}_{\chi_{\kappa}(U^{df})}.
\end{equation}

For the equations of \eqref{6.13n}, applying $\Gamma^\alpha$ on $\eqref{6.13n}_1$ and $\partial^a\Gamma^\alpha$ $(a=0,1)$ on $\eqref{6.13n}_2$, respectively, we come to
\begin{equation}
\left\{
    \begin{aligned}
    \partial_t^2 \Gamma^\alpha U^{cf}-c_1^2\Delta \Gamma^\alpha U^{cf}=&\sum_{|\beta+\gamma|\leq|\alpha|}C_\alpha^\beta  
     \nabla\Bigg[g_1\left(\nabla\cdot\Gamma^\beta U^{cf}\right)\left(\nabla^\bot\cdot \Gamma^\gamma U^{df}\right)\\
      &\qquad\qquad\qquad+g_2\left(\nabla^\bot\cdot\Gamma^\beta U^{df}\right)\left(\nabla^\bot\cdot\Gamma^\gamma U^{df}\right)\Bigg]\\
    &+\sum_{|\alpha+\beta+\delta|\leq|\alpha|}C_\alpha^\beta C_{\gamma+\delta}^\gamma\nabla\Bigg[h_1\left(\nabla\cdot\Gamma^\beta U^{cf}\right)\left(\nabla\cdot\Gamma^\gamma U^{cf}\right)\left(\nabla^\bot\cdot\Gamma^\delta U^{df}\right)\\
    &\qquad\qquad\qquad\qquad+h_2\left(\nabla\cdot \Gamma^\beta U^{cf}\right)\left(\nabla^\bot\cdot \Gamma^\gamma U^{df}\right)\left(\nabla^\bot\cdot \Gamma^\delta U^{df}\right)\\
    &\qquad\qquad\qquad\qquad+h_3\left(\nabla^\bot\cdot \Gamma^\beta U^{df}\right)\left(\nabla^\bot\cdot \Gamma^\gamma U^{df}\right)\left(\nabla^\bot\cdot \Gamma^\delta U^{df}\right)\Bigg]\\
   \partial_t^2 \partial\Gamma^\alpha U^{df}-c_2^2\Delta \partial\Gamma^\alpha  U^{df}=&0.
\end{aligned}\right.
\end{equation}

Denote the Morawetz type multipliers:
\begin{equation}
    Z_1\Gamma^\alpha U^{cf}=\partial_t\Gamma^\alpha U^{cf}+\left(t+\frac{1}{c_1}\right)\partial_t\Gamma^\alpha U^{cf}+r\partial_r\Gamma^\alpha U^{cf}+\frac12\Gamma^\alpha U^{cf},
\end{equation}
\begin{equation}
    Z_2\partial^a\Gamma^\alpha U^{df}=\partial_t\partial^a\Gamma^\alpha U^{df}+\left(t+\frac{1}{c_2}\right)\partial_t\partial^a\Gamma^\alpha U^{df}+r\partial_r\partial^a\Gamma^\alpha U^{df}+\frac12\partial^a\Gamma^\alpha U^{df}\qquad a=0,1,
\end{equation}
then the related Morawetz-type energy will be deduced:
\begin{equation}\label{n6.23}
    M_{\kappa,\kappa+1}(U^{cf},U^{df})
    =M_{\kappa}(U^{cf})+M_{\kappa+1}(U^{df}),
\end{equation}
where
\begin{equation}
\begin{aligned}
    M_{\kappa}(U^{cf})
    =&\frac12\sum_{|\alpha|\leq\kappa-1}\int_{\mathbb{R}^2}\Bigg[\left(t+\frac{1}{c_1}+1\right)|\partial_t\Gamma^\alpha U^{cf}|^2+c_1^2\left(t+\frac{1}{c_1}+1\right)|\nabla \Gamma^\alpha U^{cf}|^2\\
    &\qquad\qquad+2r\partial_r\Gamma^\alpha U^{cf}\partial_t\Gamma^\alpha U^{cf}+\Gamma^\alpha U^{cf}\partial_t\Gamma^\alpha U^{cf}\Bigg]dx,\\
        M_{\kappa+1}(U^{df})=&\frac12\sum_{|\alpha|\leq\kappa-1}\int_{\mathbb{R}^2}\Bigg[\left(t+\frac{1}{c_2}+1\right)|\partial_t\Gamma^\alpha U^{df}|^2+c_2^2\left(t+\frac{1}{c_2}+1\right)|\nabla \Gamma^\alpha U^{df}|^2\\
       &\qquad\qquad+2r\partial_r\Gamma^\alpha U^{df}\partial_t\Gamma^\alpha U^{df}+\Gamma^\alpha U^{df}\partial_t\Gamma^\alpha U^{df}\Bigg]dx\\
    &+\frac12\sum_{|\alpha|\leq\kappa-1}\int_{\mathbb{R}^2}\Bigg[\left(t+\frac{1}{c_2}+1\right)|\partial_t\partial\Gamma^\alpha U^{df}|^2+c_2^2\left(t+\frac{1}{c_2}+1\right)|\nabla \partial\Gamma^\alpha U^{df}|^2\\
    &\qquad\qquad+2r\partial_r\partial\Gamma^\alpha U^{df}\partial_t\partial\Gamma^\alpha U^{df}+\partial\Gamma^\alpha U^{df}\partial_t\partial\Gamma^\alpha U^{df}\Bigg]dx.
    \end{aligned}
\end{equation}
\begin{remark}
    Noting that we need one more derivative $\partial$ for the curl-part $U^{df}$ than for the divergence-part $U^{cf}$. This is because the highest-order derivative terms of $U^{df}$ emerges in the equation of $U^{cf}$, which means that the system of $(U^{cf},U^{df})$ is a non-symmetric hyperbolic system. This requirement is coincide with the scenarios of low-regularity local well-posedness for Euler equations \cite{Wang} and admissible harmonic elasticity \cite{ACY}.
\end{remark}

Similar to the scalar quasilinear wave equation case, we can show that the Morawetz type energy \eqref{n6.23} is non-negative, and its lower bound is demonstrated in the following lemma.
\begin{lem}\label{lem62}
   If the initial data $(U_0,U_1)$ is compactly supported in $\{x:|x|\leq 1^-\}$, then the Morawetz-type energy \eqref{n6.23} is non-negative and can be bounded from below by
    \begin{equation}\label{621}
        M_{\kappa,\kappa+1}(U^{cf},U^{df})\geq E_{\kappa,\kappa+1}(U^{cf},U^{df})+\tilde{M}_{\kappa,\kappa+1}(U^{cf},U^{df}),
    \end{equation}
    where
    \begin{equation}
        \begin{aligned}
        \tilde{M}_{\kappa,\kappa+1}(U^{cf},U^{df})=& 
        \sum_{|\alpha|\leq\kappa-1}\int_{\mathbb{R}^2}\frac{c_1t+r+1}{8c_1}\left(\partial_t\Gamma^\alpha U^{cf}+c_1\partial_r\Gamma^\alpha U^{cf}\right)^2dx\\
  &+\sum_{|\alpha|\leq\kappa-1}\int_{\mathbb{R}^2}\frac{c_1t-r+1}{8c_1}\left(\partial_t\Gamma^\alpha U^{cf}-c_1\partial_r\Gamma^\alpha U^{cf}\right)^2dx\\
   &+\sum_{|\alpha|\leq\kappa-1}\int_{\mathbb{R}^2}\frac{c_1(c_1t+1)}{2r^2}\left(\Omega \Gamma^\alpha U^{cf}\right)^2dx\\
   &+\sum_{|\alpha|\leq\kappa-1}\int_{\mathbb{R}^2}\frac{c_1}{8(c_1t+1)}|\Gamma^\alpha U^{cf}|^2dx\\
  & +\sum_{a=0,1\atop |\alpha|\leq\kappa-1}\int_{\mathbb{R}^2}\frac{c_2t+r+1}{8c_2}\left(\partial_t\partial^a\Gamma^\alpha U^{df}+c_2\partial_r\partial^a\Gamma^\alpha U^{df}\right)^2dx\\
  &+\sum_{a=0,1\atop |\alpha|\leq\kappa-1}\int_{\mathbb{R}^2}\frac{c_2t-r+1}{8c_2}\left(\partial_t\partial^a\Gamma^\alpha U^{df}-c_2\partial_r\partial^a\Gamma^\alpha U^{df}\right)^2dx\\
   &+\sum_{a=0,1\atop |\alpha|\leq\kappa-1}\int_{\mathbb{R}^2}\frac{c_2(c_2t+1)}{2r^2}\left(\Omega \partial^a\Gamma^\alpha U^{df}\right)^2dx\\
   &+\sum_{a=0,1\atop |\alpha|\leq\kappa-1}\int_{\mathbb{R}^2}\frac{c_2}{8(c_2t+1)}|\partial^a\Gamma^\alpha U^{df}|^2dx.
        \end{aligned}
    \end{equation}
\end{lem}
Moreover, the Morawetz-type energy \eqref{n6.23} is also bounded from above by $E(U)+\tilde{M}(U)$, that is
\begin{lem}
     If the initial data $(U_0,U_1)$ is compactly supported in $\{x:|x|\leq 1^-\}$, then
       \begin{equation}
        M_{\kappa,\kappa+1}(U^{cf},U^{df})\approx E_{\kappa,\kappa+1}(U^{cf},U^{df})+\tilde{M}_{\kappa,\kappa+1}(U^{cf},U^{df}).
    \end{equation}  
\end{lem}

Before showing the proof of Theorem \ref{thmhar}, we need some preliminaries.
\begin{lem}\label{lem63}
   For any 2-vector field $U\in C_c^\infty(\mathbb{R}^2)$, we have
  \begin{equation}\label{615}
     \Delta U=\nabla(\nabla\cdot U)+\nabla^\bot(\nabla^\bot\cdot U).
 \end{equation}
\end{lem}
\emph{Proof.} Noting that $\nabla^\bot=(\partial_2, -\partial_1)$, the decomposition \eqref{615} can be confirmed by calculation directly.

The following result can be referred to Lemma 4.5 in \cite{DLMS}, see also Lemma 6.4 in \cite{ACY}.
\begin{lem}
    For any 2-vector field $U\in C_c^\infty(\mathbb{R}^2)$, the following estimate holds:
    \begin{equation}
        \|\partial U\|_{L^2(\mathbb{R}^2)}\leq \|\nabla\cdot U\|_{L^2(\mathbb{R}^2)}+\|\nabla^\bot\cdot U\|_{L^2(\mathbb{R}^2)}.
    \end{equation}
\end{lem}
\emph{Proof.} The above estimate can be obtained via integrating by parts directly.

When we carry on energy estimates, the following $L^\infty-L^2$ estimate, see Lemma 2.5 in Li \cite{Ld}, will be used.
\begin{lem}\label{lem67}
    Let $U\in C_c^\infty(\mathbb{R}^2)$ be a solution of system \eqref{hweq}, we have
    \begin{equation}
    \langle t\rangle^{1/2}\|\partial U^{cf}\|_{L^\infty(\mathbb{R}^2)}\lesssim E_2^{1/2}(U^{cf})+\sum_{a=0,1}\left\|\langle c_1t-r\rangle\partial\nabla\Gamma^a U^{cf}\right\|_{L^2(\mathbb{R}^2)},
    \end{equation}
    \begin{equation}
    \langle t\rangle^{1/2}\|\partial U^{df}\|_{L^\infty(\mathbb{R}^2)}\lesssim E_2^{1/2}(U^{df})+\sum_{a=0,1}\left\|\langle c_2t-r\rangle\partial\nabla\Gamma^a U^{df}\right\|_{L^2(\mathbb{R}^2)}.
    \end{equation}
\end{lem}

\subsection{Weighted $L^2$ estimates}
The following Klainerman-Sideris' inequality for the 2D admissible harmonic equations is also true, we refer the reader to Peng-Zha \cite{PZ}.
\begin{lem}\label{lem6.7}
For any $U\in C^\infty([0,T]\times \mathbb{R}^2)$ and $\kappa\geq 2$, there holds
\begin{equation}
    \chi_{\kappa}(U)\lesssim E_{\kappa}^{1/2}(U)+t\sum_{|\alpha|\leq \kappa-2}\left\|L\Gamma^\alpha U\right\|_{L^2},
\end{equation}
where $LU=\partial_t^2 U-c_2^2\Delta U-(c_1^2-c_2^2)\nabla(\nabla\cdot U)$.
\end{lem}

Recalling that $U=U^{cf}+U^{df}$  satisfying $\nabla^\bot\cdot U^{cf}=0$ and $\nabla\cdot U^{df}=0$, by using Lemma \ref{lemM1}, then it follows 
\begin{lem}\label{L68}
    Let $U^{cf}$ and $U^{df}$ be defined as in \eqref{611}, we then have $\textup{supp} U^{cf}(\cdot,t)\subset\{x:|x|\leq c_1t+1\}$ and $\textup{supp} U^{df}(\cdot,t)\subset\{x:|x|\leq c_2t+1\}$. Moreover, 
    \begin{equation}\label{6.39}
        \left\|\langle c_1t-r\rangle^{1/2}\partial U^{cf}\right\|_{L^2}\leq M_1^{1/2}(U^{cf})
    \end{equation}
    and
    \begin{equation}\label{6.40}
        \left\|\langle c_2t-r\rangle^{1/2}\partial U^{df}\right\|_{L^2}\leq M_1^{1/2}(U^{df}). 
    \end{equation}
\end{lem}

For the 2D quasilinear wave system with multiple speeds and quadratic nonlinear terms, the weighted general energy $\chi_{\kappa}(U)$ cannot be controlled by $E_\kappa^{1/2}(U)$ directly. Fortunately, with the Morawetz type energy in hand, we have 
\begin{lem}\label{lem6.9}
    Let $U\in C_c^\infty([0,T]\times \mathbb{R}^2)$ be a solution of \eqref{hw}-\eqref{hwdata}. Let $\kappa>14$. Assume that $E_{[\kappa/2]+5}(U)\ll 1$. If $g_0=g_2=h_0=0$, i.e. the nonlinear terms of \eqref{hw} are composed of some quadratic and cubic terms, like 
    \[
    \nabla\left[g_1\left(\nabla\cdot U\right)\left(\nabla^\bot\cdot U\right)\right]+\nabla\left[h_1\left(\nabla\cdot U\right)^2\left(\nabla^\bot\cdot U\right)+h_2\left(\nabla\cdot U\right)\left(\nabla^\bot\cdot U\right)^2+h_3\left(\nabla^\bot\cdot U\right)^3\right], 
    \]
    then it holds
    \begin{equation}\label{n6.35}
    \begin{aligned}
    \chi_\kappa(U)\lesssim E_{\kappa}^{1/2}(U)+M_{\kappa}^{1/2}(U^{df})\left[E_{[\kappa/2]+5}^{1/2}(U)+\chi_{[\kappa/2]+5}(U)\right]+E_{[\kappa/2]+5}^{1/2}(U^{df})M_\kappa^{1/2}(U^{cf}).
        \end{aligned}\end{equation}
       In particular, for the weighed $L^2$ norm of the divergence-free part, it holds
       \begin{equation}\label{n6.36}
        \chi_\kappa(U^{df})\lesssim E_{\kappa}^{1/2}(U^{df}).
       \end{equation}
\end{lem}
\emph{Proof.} Since $U^{df}$ obeys a single decoupled linear wave equation, the null condition is satisfied automatically and \eqref{n6.36} holds obviously by Klainerman-Sideris' inequality in Lemma \ref{lem6.7}. For \eqref{n6.35}, first we have
\begin{equation}\label{6.34}
\begin{split}
    \chi_\kappa(U^{cf})\lesssim& E_\kappa^{1/2}(U^{cf})+ \sum_{|\alpha|\leq\kappa-2}t\left\|\Box_{c_1}\Gamma^\alpha U^{cf}\right\|_{L^2}.
\end{split}\end{equation}
By \eqref{6.39}, the quadratic nonlinear terms in \eqref{6.34} can be estimated as
\begin{equation}\label{6.35}
\begin{split}
t\left\|\Box_{c_1}\Gamma^\alpha U^{cf}\right\|_{L^2}\lesssim&  t\left\|g_1\nabla^2\Gamma^\beta U^{cf}\nabla^\bot\cdot \Gamma^\gamma U^{df}\right\|_{L^2}+ t\left\|g_1\nabla\cdot\Gamma^\beta U^{cf}\nabla^2 \Gamma^\gamma U^{df}\right\|_{L^2}\\
\lesssim&\left\{\begin{aligned}
&\left\|\langle c_1t-r\rangle\nabla^2\Gamma^\beta U^{cf}\right\|_{L^2}\left\|\langle r\rangle^{1/2}\langle c_2t-r\rangle^{1/2}\nabla^\bot\cdot\Gamma^\gamma U^{df}\right\|_{L^\infty}\\
&+\left\|\langle c_1t-r\rangle^{1/2}\nabla\cdot\Gamma^\beta U^{cf}\right\|_{L^2}\left\|\langle r\rangle^{1/2}\langle c_2t-r\rangle^{1/2}\nabla^2\Gamma^\gamma U^{df}\right\|_{L^\infty}, \ \text{if}\ |\beta|\geq|\gamma|,\\
&\left\|\langle r\rangle^{1/2}\langle c_1t-r\rangle^{1/2}\nabla^2\Gamma^\beta U^{cf}\right\|_{L^\infty}\left\|\langle c_2t-r\rangle^{1/2}\nabla^\bot\cdot\Gamma^\gamma U^{df}\right\|_{L^2}\\
&+\left\|\langle r\rangle^{1/2}\langle c_1t-r\rangle^{1/2}\nabla\cdot\Gamma^\beta U^{cf}\right\|_{L^\infty}\left\|\langle c_2t-r\rangle\nabla^2\Gamma^\gamma U^{df}\right\|_{L^2}, \ \text{if}\ |\beta|\leq|\gamma|
\end{aligned}\right.\\
\lesssim&E_{[\kappa/2]+4}^{1/2}(U^{df})\chi_\kappa(U^{cf})+E_{[\kappa/2]+5}^{1/2}(U^{df})M_\kappa^{1/2}(U^{cf})\\
&+\left(E_{[\kappa/2]+5}^{1/2}(U^{cf})+\chi_{[\kappa/2]+5}(U^{cf})\right)\left(M_\kappa^{1/2}(U^{df})+E_\kappa^{1/2}(U^{df})\right).
\end{split}\end{equation}
For the cubic nonlinear terms, we only need to analyze the term $\nabla(\nabla^\bot\cdot U^{df})^3$, and the other cubic terms can be estimated in a similar way as used for the preceding quadratic nonlinear terms, with an additional factor $E_{[\kappa/3]+4}^{1/2}(U)$ arising. Using Lemma \ref{lem67}, we have
\begin{equation}\label{636}
\begin{aligned}
&t\left\|h_3\nabla\left[\left(\nabla^\bot\cdot \Gamma^\beta U^{df}\right)\left(\nabla^\bot\cdot \Gamma^\gamma U^{df}\right)\left(\nabla^\bot\cdot \Gamma^\delta U^{df}\right)\right]\right\|_{L^2}\\
\lesssim& \left\|t^{1/2}\nabla^\bot\cdot \Gamma^{[\kappa/3]+1} U^{df}\right\|_{L^\infty} \left\|t^{1/2}\nabla^\bot\cdot \Gamma^{[\kappa/3]+1} U^{df}\right\|_{L^\infty} \left\|\nabla^\bot\cdot \Gamma^{\kappa-1} U^{df}\right\|_{L^2}\\
\lesssim&\left(E_{[\kappa/3]+4}^{1/2}(U^{df})+\chi_{[\kappa/3]+4}(U^{df})\right)^2E_\kappa^{1/2}(U^{df})\\
\lesssim&E_{[\kappa/3]+4}(U^{df})E_\kappa^{1/2}(U^{df}).
\end{aligned}
\end{equation}
Combing with \eqref{6.34},  \eqref{6.35} and \eqref{636}, we complete the proof of \eqref{n6.35} if $E_{[\kappa/2]+5}^{1/2}$ is small enough.\hfill$\Box$

In the above lemmas, we do not use the Lorentz boost. The Morawetz type energy estimate is still valid when the Lorentz boost is absent. But in the process of estimating the classical general energy, we employ the Lorentz boost only for the equation of $U^{df}$. That is, if we apply $\Gamma$ on $\Box_{c_2}U^{df}=0$, this $\Gamma$ includes the vector fields $L_i=c_2t\partial_i+\frac{x_i}{c_2}\partial_t$ ($i=1,2$).  Without misunderstanding, we also use notation $\Gamma$ for $\Gamma U^{df}$.
\begin{lem}\label{lem611}
    It holds
    \begin{equation}
        \left|(c_2t-r)\partial U^{df}\right|\leq \left|\Gamma U^{df}\right|.
    \end{equation}
\end{lem}
\emph{Proof.} One can find the original result in Klainerman \cite{Kl85}. The main step of the proof relies on the observation 
\begin{equation}
    \left(c_2^2t^2-r^2\right)\partial_r=c_2tL_r-rS,
\end{equation}
where $L_r=\sum_{i=1}^2w_iL_i=c_2t\partial_r+\frac{r}{c_2}\partial_t$.\hfill$\Box$

For the divergence-free part $U^{df}$, the following $L^2$ estimate (see Lemma 4.3.1 in Li-Zhou \cite{LZ}) also plays a crucial role in controlling the solution itself.
\begin{lem}\label{L2L}
    Let $U^{df}$ be the solution of $\eqref{6.13n}_2$ with initial data $$U^{df}(0,x)=U^{df}_0(x),\ U^{df}_t(0,x)=U^{df}_1(x).$$ Then we have
    \begin{equation}\label{L2L2}
\left\|U^{df}\right\|_{L^2(\mathbb{R}^2)}\leq\left\|U_0^{df}\right\|_{L^2(\mathbb{R}^2)}+C\sqrt{\ln(3+t)}\left\|\langle r\rangle^2 U_1^{df}\right\|_{L^2(\mathbb{R}^2)}.
    \end{equation}
\end{lem}

\subsection{Energy estimates}
Let $\kappa\geq 14$ . For the admissible harmonic elastodynamics, we shall build the related Morawetz type energy estimate and classical general energy a prior estimate:
\begin{equation}
    \frac{d}{dt}M_{\kappa,\kappa+1}(U^{cf},U^{df})\lesssim \langle t\rangle^{-1}E_{\mu,\mu+1}^{1/2}(U^{cf},U^{df})M_{\kappa,\kappa+1}(U^{cf},U^{df}),\qquad 0<t<T,
\end{equation}
and 
\begin{equation}
    \frac{d}{dt}E_{\mu,\mu+1}(U^{cf},U^{df})\lesssim \langle t\rangle^{-3/2}\ln^{1/2}(3+t)E_{\mu,\mu+1}(U^{cf},U^{df})M_{\kappa,\kappa+1}^{1/2}(U^{cf},U^{df}),\qquad 0<t<T.
\end{equation}

{\bf Morawetz-type energy estimate:} Let $\kappa\geq 14$ and $\mu=\kappa-2$. For $|\alpha|\leq\kappa-1$, multiplying on the both sides of $\eqref{6.13n}_2$ with $Z_2\partial^a\Gamma^\alpha U^{df}$ $(a=0,1)$ yields
\begin{equation}\label{6.42}
 \frac{d}{dt}M_{\kappa+1}(U^{df})=0.
\end{equation}

 Multiplying the both sides of $\eqref{6.13n}_1$ with $Z_1\Gamma^\alpha U^{cf}$, and making integration by parts, we get:
\begin{equation}\label{6.43}
    \begin{split}
        \frac{d}{dt}\tilde{M}_{\kappa}(U^{cf})\simeq&\underbrace{\frac12g_1\sum_{|\alpha|=0}^{\kappa-1}\int_{\mathbb{R}^2}\left(\nabla^\bot\cdot \Gamma U^{df}+\nabla^\bot\cdot  U^{df}+\left(1+\frac{1}{c_1}\right)\partial_t\nabla^\bot\cdot U^{df}\right)\left|\nabla\cdot\Gamma^\alpha U^{cf}\right|^2 dx}_{J_1}\\
        &+\underbrace{g_1\sum_{|\alpha|=0}^{\kappa-1}\int_{\mathbb{R}^2}\nabla\left(\nabla\cdot U^{cf} \nabla^\bot\cdot\Gamma^\alpha U^{df}\right)\partial_t\Gamma^\alpha U^{cf}dx}_{J_2}\\
        &+\underbrace{g_1\sum_{|\alpha|=0}^{\kappa-1}\int_{\mathbb{R}^2}\nabla\cdot U^{cf} \left(2\nabla^\bot\cdot\Gamma^\alpha U^{df}+\nabla^\bot\cdot\Gamma^{\alpha+1} U^{df}\right)\nabla\cdot\Gamma^\alpha U^{cf}dx}_{J_3}\\
        &+ \underbrace{g_1\sum_{|\alpha|=0}^{\kappa-1}\sum_{|\beta+\gamma|\leq|\alpha|\atop |\beta|,|\gamma|<|\alpha|}C_\alpha^\beta \int_{\mathbb{R}^2}\nabla\left(\nabla\cdot \Gamma^\beta U^{cf} \nabla^\bot\cdot\Gamma^\gamma U^{df}\right)\partial_t\Gamma^\alpha U^{cf}dx}_{J_4}\\
        &+ \underbrace{ g_1\sum_{|\alpha|=0}^{\kappa-1}\sum_{|\beta+\gamma|\leq|\alpha|\atop |\beta|,|\gamma|<|\alpha|}C_\alpha^\beta\int_{\mathbb{R}^2}\nabla\cdot \Gamma^\beta U^{cf}\left (2\nabla^\bot\cdot\Gamma^\gamma U^{df}+\nabla^\bot\cdot\Gamma^{\gamma+1} U^{df}\right)\nabla\cdot\Gamma^\alpha U^{cf}dx}_{J_5}\\
        &+ \underbrace{g_1\sum_{|\alpha|=0}^{\kappa-1}\sum_{|\beta+\gamma|\leq|\alpha|\atop |\beta|,|\gamma|<|\alpha|}C_\alpha^\beta \int_{\mathbb{R}^2}\nabla\cdot \Gamma^{\beta+1} U^{cf} \nabla^\bot\cdot\Gamma^\gamma U^{df}\nabla\cdot\Gamma^\alpha U^{cf}dx}_{J_6}\\
        &+\underbrace{g_2\sum_{|\alpha|=0}^{\kappa-1}\sum_{|\beta+\gamma|\leq|\alpha|}C_\alpha^\beta\int_{\mathbb{R}^2} \nabla^\bot\cdot \Gamma^{\beta+1} U^{df} \nabla^\bot\cdot\Gamma^\gamma U^{df}\nabla\cdot\Gamma^\alpha U^{cf}dx}_{J_7}+cubic\ terms,
    \end{split}
\end{equation}
where
\begin{equation}
    \begin{split}
       \tilde{M}_{\kappa}(U^{cf})=&M_{\kappa}(U^{cf})+\underbrace{\frac12g_1\sum_{|\alpha|=0}^{\kappa-1}\int_{\mathbb{R}^2}\left(t+1+\frac{1}{c_1}\right)\nabla^\bot\cdot U^{df}\left|\nabla\cdot\Gamma^\alpha U^{cf}\right|^2 dx}_{I_1}\\
        &+g_1\underbrace{\sum_{|\alpha|=0}^{\kappa-1}\int_{\mathbb{R}^2}\left(t+1+\frac{1}{c_1}\right) \nabla\cdot U^{cf} \nabla^\bot\cdot\Gamma^\alpha U^{df}\nabla\cdot\Gamma^\alpha U^{cf}dx}_{I_2}\\
        &+ g_1\underbrace{\sum_{|\alpha|=0}^{\kappa-1}\sum_{|\beta+\gamma|\leq|\alpha|\atop |\beta|,|\gamma|<|\alpha|}C_\alpha^\beta\int_{\mathbb{R}^2}\left(t+1+\frac{1}{c_1}\right) \nabla\cdot \Gamma^\beta U^{cf} \nabla^\bot\cdot\Gamma^\gamma U^{df}\nabla\cdot\Gamma^\alpha U^{cf}dx}_{I_3}\\
   &+g_2\underbrace{\sum_{|\alpha|=0}^{\kappa-1}\sum_{|\beta+\gamma|\leq|\alpha|}C_\alpha^\beta\int_{\mathbb{R}^2}\left(t+1+\frac{1}{c_1}\right) \nabla^\bot\cdot \Gamma^\beta U^{df} \nabla^\bot\cdot\Gamma^\gamma U^{df}\nabla\cdot\Gamma^\alpha U^{cf}dx}_{I_4}.
    \end{split}
\end{equation}
In this part, we only use $\chi_\mu(U^{df})$ without $\chi_\mu(U^{cf})$, which is bounded by $E_\mu^{1/2}(U^{df})$ in Lemma \ref{lem6.9}. In detail,
since $$\langle t\rangle\lesssim \langle r\rangle^{1/2}\langle c_1t-r\rangle^{1/2}\langle c_2t-r\rangle^{1/2},$$
we may estimate $I_1,I_2$ terms as 
\begin{equation}
\begin{split}
I_1\lesssim& \left\|\langle r\rangle^{1/2}\langle c_2t-r\rangle^{1/2}\nabla^\bot\cdot U^{df}\right\|_{L^\infty}\left\|\langle c_1t-r\rangle^{1/2}\nabla\cdot\Gamma^\alpha U^{cf}\right\|_{L^2}\left\|\nabla\cdot\Gamma^\alpha U^{cf}\right\|_{L^2}\\
\lesssim&\left(E_\mu^{1/2}(U^{df})+\chi_{\mu}(U^{df})\right)M_{\kappa}^{1/2}(U^{cf})E_{\kappa}^{1/2}(U^{cf}),
\end{split}
\end{equation}
\begin{equation}
\begin{split}
I_2\lesssim& \left\|\langle r\rangle^{1/2}\nabla\cdot U^{cf}\right\|_{L^\infty}\left\|\langle c_2t-r\rangle^{1/2}\nabla^\bot\cdot \Gamma^\alpha U^{df}\right\|_{L^2}\left\|\langle c_1t-r\rangle^{1/2}\nabla\cdot\Gamma^\alpha U^{cf}\right\|_{L^2}\\
\lesssim&E_\mu^{1/2}(U^{cf})M_{\kappa,\kappa}^{1/2}(U^{cf},U^{df}).
\end{split}
\end{equation}
For $I_3$, if $|\beta|<|\gamma|$, it holds
\begin{equation}
\begin{split}
I_3\lesssim& \left\|\langle r\rangle^{1/2}\nabla\cdot\Gamma^\beta U^{cf}\right\|_{L^\infty}\left\|\langle c_2t-r\rangle^{1/2}\nabla^\bot\cdot \Gamma^\gamma U^{df}\right\|_{L^2}\left\|\langle c_1t-r\rangle^{1/2}\nabla\cdot\Gamma^\alpha U^{cf}\right\|_{L^2}\\
\lesssim&E_\mu^{1/2}(U^{cf})M_{\kappa,\kappa}^{1/2}(U^{cf},U^{df}),
\end{split}
\end{equation}
while if $|\beta|>|\gamma|$ in $I_3$, we have
\begin{equation}
\begin{split}
I_3\lesssim& \left\|\nabla\cdot\Gamma^\beta U^{cf}\right\|_{L^2}\left\|\langle r\rangle^{1/2}\langle c_2t-r\rangle^{1/2}\nabla^\bot\cdot \Gamma^\gamma U^{df}\right\|_{L^\infty}\left\|\langle c_1t-r\rangle^{1/2}\nabla\cdot\Gamma^\alpha U^{cf}\right\|_{L^2}\\
\lesssim&\left(E_\mu^{1/2}(U^{df})+\chi_{\mu}(U^{df})\right)M_{\kappa}^{1/2}(U^{cf})E_{\kappa}^{1/2}(U^{cf}).
\end{split}
\end{equation}
Similarly, for $I_4$ it holds
\begin{equation}
\begin{split}
I_4\lesssim&\left\|\langle r\rangle^{1/2}\langle c_2t-r\rangle^{1/2}\nabla^\bot\cdot \Gamma^{[\alpha/2]+1} U^{df}\right\|_{L^\infty} \left\|\nabla\cdot\Gamma^\alpha U^{df}\right\|_{L^2}\left\|\langle c_1t-r\rangle^{1/2}\nabla\cdot\Gamma^\alpha U^{cf}\right\|_{L^2}\\
\lesssim&\left(E_\mu^{1/2}(U^{df})+\chi_{\mu}(U^{df})\right)M_{\kappa}^{1/2}(U^{cf})E_{\kappa}^{1/2}(U^{df}).
\end{split}
\end{equation}
Thanks to \eqref{n6.36}, we have $\chi_{\mu}(U^{df})\lesssim E_\mu^{1/2}(U^{df})$. Finally, the terms $I_1-I_4$ are bounded by
\begin{equation}
I_1, I_2,I_3,I_4\lesssim E_\mu^{1/2}(U)M_{\kappa,\kappa}(U^{cf}, U^{df}).
\end{equation}
Adding \eqref{6.42} and \eqref{6.43} together, it is easy to see $I_1-I_3$ are small perturbations of $M_{\kappa,\kappa}(U^{cf}, U^{df})$, provided that $E_\mu^{1/2}(U)\ll \epsilon_0$ is small enough, thus 
\begin{equation}\label{6.49}
        \frac{d}{dt}\tilde{M}_{\kappa,\kappa+1}(U^{cf},U^{df})
   \simeq J_1+J_2+\cdots+J_7+cubic\ terms,
\end{equation}
where $$\tilde{M}_{\kappa,\kappa+1}(U^{cf},U^{df})\simeq M_{\kappa,\kappa+1}(U^{cf},U^{df}).$$

We are going to estimate $J_1-J_7$ on the right hand side of \eqref{6.49}. Similar to $I_1$, $J_1$ is bounded by
\begin{equation}
J_1\lesssim \frac{1}{\langle t\rangle}\left[E_\mu^{1/2}(U^{df})+\chi_{\mu}(U^{df})\right]M_{\kappa}^{1/2}(U^{cf})E_{\kappa}^{1/2}(U^{cf})\lesssim  \frac{1}{\langle t\rangle}E_\mu^{1/2}(U)M_{\kappa,\kappa}(U^{cf}, U^{df}).
\end{equation}
The estimations of $J_2$ and $J_3$ are similar to $I_2$ and we have
\begin{equation}
J_2,J_3\lesssim  \frac{1}{\langle t\rangle}E_\mu^{1/2}(U)M_{\kappa,\kappa+1}(U^{cf}, U^{df}).
\end{equation}
And $J_4,J_5,J_6$ can be controlled as $I_3$ to get
\begin{equation}\label{6.52}
J_4,J_5,J_6\lesssim  \frac{1}{\langle t\rangle}E_\mu^{1/2}(U)M_{\kappa,\kappa}(U^{cf}, U^{df}).
\end{equation}
Similar to $I_4$, we have
\begin{equation}
J_7\lesssim  \frac{1}{\langle t\rangle}E_\mu^{1/2}(U)M_{\kappa,\kappa+1}(U^{cf}, U^{df}).
\end{equation}
The  estimations for the cubic terms are analogous to those for the quadratic terms, with an additional $L^\infty$ norm of a lower-order derivative term, which can be controlled by $E_\mu^{1/2}(U)$ via the standard Sobolev embedding inequality. From \eqref{6.49} to \eqref{6.52}, we finally obtain the higer-order Morawetz type energy estimate:
\begin{equation}\label{656}
        \frac{d}{dt}\tilde{M}_{\kappa,\kappa+1}(U^{cf},U^{df})\lesssim \frac{1}{\langle t\rangle}\left[E_\mu^{1/2}(U)+E_\mu(U)\right]M_{\kappa,\kappa+1}(U^{cf}, U^{df}).
\end{equation}

{\bf Classical general energy estimate:} In this part, we consider that the admissible harmonic elastic waves \eqref{hweq} satisfying $g_0=g_2=h_0=0$. The dynamic of the curl-free part obeys
\begin{equation}
\begin{aligned}
  \partial_t^2 \Gamma^\alpha U^{cf}-c_1^2\Delta \Gamma^\alpha U^{cf}=&\sum_{|\beta+\gamma|\leq|\alpha|}C_\alpha^\beta  
    g_1 \nabla\left[\left(\nabla\cdot\Gamma^\beta U^{cf}\right)\left(\nabla^\bot\cdot \Gamma^\gamma U^{df}\right)\right]\\
&+\sum_{|\beta+\gamma+\delta|\leq|\alpha|}C_\alpha^\beta C_{\gamma+\delta}^\gamma\nabla\Bigg[h_1\left(\nabla\cdot\Gamma^\beta U^{cf}\right)\left(\nabla\cdot\Gamma^\gamma U^{cf}\right)\left(\nabla^\bot\cdot\Gamma^\delta U^{df}\right)\\
    &\qquad\qquad\qquad\qquad+h_2\left(\nabla\cdot \Gamma^\beta U^{cf}\right)\left(\nabla^\bot\cdot \Gamma^\gamma U^{df}\right)\left(\nabla^\bot\cdot \Gamma^\delta U^{df}\right)\\
    &\qquad\qquad\qquad\qquad+h_3\left(\nabla^\bot\cdot \Gamma^\beta U^{df}\right)\left(\nabla^\bot\cdot \Gamma^\gamma U^{df}\right)\left(\nabla^\bot\cdot \Gamma^\delta U^{df}\right)\Bigg].
\end{aligned}
\end{equation}

Let $\mu\leq\kappa-2$ and $|\alpha|\leq\mu-1$, multiplying the both sides of $\eqref{hweq}_2$ with $\partial_t\Gamma^{\alpha+1} U^{df}$, we get
\begin{equation}\label{6.53}
 \frac{d}{dt}E_{\mu+1}(U^{df})=0,
\end{equation}
while multiplying the both sides of $\eqref{hweq}_1$ with $\partial_t\Gamma^\alpha U^{cf}$ we get
\begin{equation}
\begin{aligned}
  \frac{d \tilde{E}_{\mu}(U^{cf})}{dt}=&\underbrace{\frac12\sum_{|\alpha|\leq\mu-1}\left\|g_1\left(\nabla^\bot\cdot U^{df}\right)\left(\nabla\nabla\cdot\Gamma^\alpha U^{cf}\right)\left(\partial_t\Gamma^\alpha U^{cf}\right)\right\|_{L^1}}_{H_1}\\
  &+\underbrace{\sum_{|\alpha|\leq\mu-1}\sum_{|\beta+\gamma|\leq|\alpha|\atop |\beta|<|\alpha|}C_\alpha^\beta \left\|g_1\partial_t\left(\nabla\cdot\Gamma^\beta U^{cf}\nabla^\bot\cdot\Gamma^\gamma U^{df}\right)\nabla\cdot\Gamma^\alpha U^{cf}\right\|_{L^1}}_{H_2}\\
  &+\underbrace{\frac12\sum_{|\alpha|\leq\mu-1}\left\|h_1\partial_t\left(\nabla\cdot U^{cf}\nabla^\bot\cdot U^{df}\right)\left(\nabla\cdot\Gamma^\alpha U^{cf}\right)^2\right\|_{L^1}}_{N_1}\\
  &+\underbrace{\sum_{|\alpha|\leq\mu-1}\left\|h_2\nabla^\bot\cdot U^{df}\partial_t \nabla^\bot\cdot U^{df}\left(\nabla\cdot\Gamma^\alpha U^{cf}\right)^2\right\|_{L^1}}_{N_2}\\
  &+ \underbrace{\sum_{|\alpha|\leq\mu-1}\sum_{|\beta+\gamma+\delta|\leq|\alpha|\atop |\beta|<|\alpha|}C_\alpha^\beta C_{\gamma+\delta}^\gamma\left\|h_1\partial_t\left(\nabla\cdot\Gamma^\beta U^{cf}\nabla\cdot\Gamma^\gamma U^{cf}\nabla^\bot\cdot\Gamma^\delta U^{df}\right)\nabla\cdot\Gamma^\alpha U^{cf}\right\|_{L^1}}_{N_3}\\
  &+ \underbrace{\sum_{|\alpha|\leq\mu-1}\sum_{|\beta+\gamma+\delta|\leq|\alpha|\atop |\beta|<|\alpha|}C_\alpha^\beta C_{\gamma+\delta}^\gamma\left\|h_2\partial_t\left(\nabla\cdot\Gamma^\beta U^{cf}\nabla^\bot\cdot\Gamma^\gamma U^{df}\nabla^\bot\cdot\Gamma^\delta U^{df}\right)\nabla\cdot\Gamma^\alpha U^{cf}\right\|_{L^1}}_{N_4},\\
\end{aligned}
\end{equation}
where 
\begin{equation}
\begin{aligned}
  \tilde{E}_{\mu}(U^{cf})=&{E}_{\mu}(U^{cf})
  +\underbrace{\sum_{|\alpha|\leq\mu-1}\sum_{|\beta+\gamma|\leq|\alpha|\atop |\beta|<|\alpha|}C_\alpha^\beta \left\|g_1\nabla^\bot\cdot\Gamma^\beta U^{df}\nabla^\bot\cdot\Gamma^\gamma U^{df}\nabla\cdot\Gamma^\alpha U^{cf}\right\|_{L^1}}_{K_1}\\
  &+\underbrace{\frac12\sum_{|\alpha|\leq\mu-1}\left\|h_1\nabla\cdot U^{cf}\nabla^\bot\cdot U^{df}\left(\nabla\cdot\Gamma^\alpha U^{cf}\right)^2\right\|_{L^1}}_{K_2}\\
  &+\underbrace{\frac12\sum_{|\alpha|\leq\mu-1}\left\|h_2\left(\nabla^\bot\cdot U^{df}\right)^2\left(\nabla\cdot\Gamma^\alpha U^{cf}\right)^2\right\|_{L^1}}_{K_3}\\
  &+ \underbrace{\sum_{|\alpha|\leq\mu-1}\sum_{|\beta+\gamma+\delta|\leq|\alpha|\atop |\beta|<|\alpha|}C_\alpha^\beta C_{\gamma+\delta}^\gamma\left\|h_1\nabla\cdot\Gamma^\beta U^{cf}\nabla\cdot\Gamma^\gamma U^{cf}\nabla^\bot\cdot\Gamma^\delta U^{df}\nabla\cdot\Gamma^\alpha U^{cf}\right\|_{L^1}}_{K_4}\\
  &+ \underbrace{\sum_{|\alpha|\leq\mu-1}\sum_{|\beta+\gamma+\delta|\leq|\alpha|\atop |\beta|<|\alpha|}C_\alpha^\beta C_{\gamma+\delta}^\gamma\left\|h_2\nabla\cdot\Gamma^\beta U^{cf}\nabla^\bot\cdot\Gamma^\gamma U^{df}\nabla^\bot\cdot\Gamma^\delta U^{df}\nabla\cdot\Gamma^\alpha U^{cf}\right\|_{L^1}}_{K_5}.
\end{aligned}
\end{equation}
It is easy to verify that
\begin{equation}
K_1,\cdots,K_5\lesssim\left(E_{[\mu/2]+4}(U^{cf})+E_{[\mu/2]+4}(U^{df})\right)\left(E_{\mu}(U^{cf})+E_{\mu}(U^{df})\right).
\end{equation}
Assume that $E_{[\mu/2]+4}(U^{cf})+E_{[\mu/2]+4}(U^{df})$ is small enough, then $\tilde{E}_{\mu,\mu+1}(U^{cf},U^{df})$ is a small perturbation of $E_{\mu,\mu+1}(U^{cf}, U^{df})$.
Together with \eqref{6.53}, we come to
\begin{equation}\label{6.57}
  \frac{d \tilde{E}_{\mu,\mu+1}(U^{cf},U^{df})}{dt}=H_1+H_2+N_1+N_2+N_3+N_4,
\end{equation}
in which 
\begin{equation}
\tilde{E}_{\mu,\mu+1}(U^{cf},U^{df})=\tilde{E}_{\mu}(U^{cf})+E_{\mu+1}(U^{df}).
\end{equation}

For $H_1$ and $H_2$, we use the fact 
\begin{equation}
    \langle t\rangle^{3/2}\lesssim \langle r\rangle^{1/2}\langle c_1t-r\rangle\langle c_2t-r\rangle
\end{equation}
to obtain 
\begin{equation}\label{H1}
\begin{aligned}
    H_1\lesssim&\frac{1}{\langle t\rangle^{3/2}}\left\|\langle r\rangle^{1/2}\langle c_2t-r\rangle\nabla^\bot\cdot U^{df}\right\|_{L^\infty}\left\|\langle c_1t-r\rangle\nabla\nabla\cdot\Gamma^\alpha U^{cf}\right\|_{L^2}\left\|\partial_t\Gamma^\alpha U^{cf}\right\|_{L^2}.
\end{aligned}\end{equation}
By \eqref{L2L2} in Lemma \ref{L2L}, it holds
\begin{equation}\label{6.65}
    \begin{aligned}
       &\left\|\langle r\rangle^{1/2}\langle c_2t-r\rangle\nabla^\bot\cdot  U^{df}\right\|_{L^\infty}\\
       \lesssim&\left\|\langle r\rangle^{1/2} \Gamma U^{df}\right\|_{L^\infty}\\
       \lesssim&E_{3}^{1/2}(U^{df})+\left\|\Gamma^{3}U^{df}\right\|_{L^2}\\
       \lesssim&E_{3}^{1/2}(U^{df})+\left\| \Gamma^{3}U^{df}(0,\cdot)\right\|_{L^2}+\ln^{1/2}(3+t)\left\|\langle r\rangle^2 \Gamma^{3}\partial_tU^{df}(0,\cdot)\right\|_{L^2}.
    \end{aligned}
\end{equation}
Plugging \eqref{6.65} into \eqref{H1} yields
\begin{equation}\label{665}
\begin{aligned}
H_1\lesssim\frac{1}{\langle t\rangle^{3/2}}\chi_\kappa(U^{cf})E_\mu^{1/2}(U^{cf})\left(E_{3}^{1/2}(U^{df})+\ln^{1/2}(3+t)\|\Lambda U^{df}(0)\|_4\right),
\end{aligned}\end{equation}
where
\begin{equation}
    \left\|\Lambda U^{df}(0)\right\|_{\mu+1}=\sum_{\gamma=0}^{\mu}\left(\left\|\Gamma^\gamma U_0^{df}\right\|_{L^2}+\left\|\langle r\rangle^2 \Gamma^{\gamma}U^{df}_1\right\|_{L^2}\right).
\end{equation}

For $H_2$, we fisrt have
\begin{equation}
\begin{aligned}
H_2\lesssim&\sum_{|\beta+\gamma|\leq|\alpha|\atop |\beta|<|\alpha|\leq\mu-1}\left\|\partial_t\nabla\cdot\Gamma^\beta U^{cf}\nabla^\bot\cdot\Gamma^\gamma U^{df}\nabla\cdot\Gamma^\alpha U^{cf}\right\|_{L^1}\\
&+\sum_{|\beta+\gamma|\leq|\alpha|\leq\mu-1\atop |\beta|<|\alpha|-1}\left\|\nabla\cdot\Gamma^\beta U^{cf}\partial_t\nabla^\bot\cdot\Gamma^\gamma U^{df}\nabla\cdot\Gamma^\alpha U^{cf}\right\|_{L^1}\\
&+\sum_{|\alpha|\leq\mu-1}\left\|\left(\partial_t\nabla^\bot\cdot\Gamma U^{df}\right)\Gamma\left( \nabla\cdot\Gamma^{\alpha-1} U^{cf}\right)^2\right\|_{L^1}, 
\end{aligned}\end{equation}
 the first two terms of $H_2$ will be discussed in two cases: $|\beta|\leq|\gamma|$ and $|\beta|>|\gamma|$, while the last term will be dealt with by integration by parts, and we come to
\begin{equation}
\begin{aligned}
H_2\lesssim&\frac{1}{\langle t\rangle^{3/2}}\Bigg(\sum_{|\beta+\gamma|\leq|\alpha|\atop |\beta|\leq|\gamma|\leq\mu-1}\left\|\langle r\rangle^{1/2}\langle c_1t-r\rangle\partial_t\nabla\cdot\Gamma^\beta U^{cf}\right\|_{L^\infty}\left\|\langle c_2t-r\rangle\nabla^\bot\cdot \Gamma^\gamma U^{df}\right\|_{L^2}\left\|\nabla\cdot\Gamma^\alpha U^{cf}\right\|_{L^2}\\
&+\sum_{|\beta+\gamma|\leq|\alpha|\atop |\gamma|<|\beta|\leq\mu-2}\left\|\langle c_1t-r\rangle\partial_t\nabla\cdot\Gamma^\beta U^{cf}\right\|_{L^2}\left\|\langle r\rangle^{1/2}\langle c_2t-r\rangle\nabla^\bot\cdot \Gamma^\gamma U^{df}\right\|_{L^\infty}\left\|\nabla\cdot\Gamma^\alpha U^{cf}\right\|_{L^2}\\
    &+\sum_{|\beta+\gamma|\leq|\alpha|\leq\mu-1\atop |\beta|<|\alpha|-1}\left\|\langle r\rangle^{1/2}\langle c_1t-r\rangle^{1/2}\nabla\cdot\Gamma^\beta U^{cf}\right\|_{L^\infty}\left\|\langle c_2t-r\rangle\partial_t\nabla^\bot\cdot \Gamma^\gamma U^{df}\right\|_{L^2}\left\|\langle c_1t-r\rangle^{1/2}\nabla\cdot\Gamma^\alpha U^{cf}\right\|_{L^2}\Bigg)\\
    &+\sum_{|\alpha|\leq\mu-1}\left\|\left(\partial_t\Gamma^2 U^{df}\right)\partial\left( \nabla\cdot\Gamma^{\alpha-1} U^{cf}\right)^2\right\|_{L^1}\\
    \lesssim&\frac{1}{\langle t\rangle^{3/2}}\chi_{[\mu/2]+4}(U^{cf})\left\| \Gamma^{\gamma+1}U^{df}\right\|_{L^2}E_\mu^{1/2}(U^{cf})\\
    &+\frac{1}{\langle t\rangle^{3/2}}\chi_{\mu}(U^{cf})\sum_{|\gamma|<[\mu/2]+1}\left\|\langle r\rangle^{1/2}\langle c_2t-r\rangle\nabla^\bot\cdot \Gamma^\gamma U^{df}\right\|_{L^\infty}E_\mu^{1/2}(U^{cf})\\
    &+\frac{1}{\langle t\rangle^{3/2}}\left(E_\mu^{1/2}(U^{cf})+\chi_\mu(U^{cf})\right)E_\mu^{1/2}(U^{df})M_\mu^{1/2}(U^{cf})\\
    &+\frac{1}{\langle t\rangle^{3/2}}\left\|\langle r\rangle^{1/2}\langle c_2t-r\rangle\partial_t \Gamma^2 U^{df}\right\|_{L^\infty}E_\mu^{1/2}(U^{cf})\chi_\mu(U^{cf}).
\end{aligned}\end{equation}
Again by \eqref{L2L2} in Lemma \ref{L2L}, it holds
\begin{equation}
    \begin{aligned}
       &\left\|\langle r\rangle^{1/2}\langle c_2t-r\rangle\nabla^\bot\cdot \Gamma^\gamma U^{df}\right\|_{L^\infty}\\
       \lesssim&\left\|\langle r\rangle^{1/2} \Gamma^{\gamma+1}U^{df}\right\|_{L^\infty}\\
       \lesssim&E_{\gamma+3}^{1/2}(U^{df})+\left\|\Gamma^{\gamma+3}U^{df}\right\|_{L^2}\\
       \lesssim&E_{\gamma+3}^{1/2}(U^{df})+\left\| \Gamma^{\gamma+3}U^{df}(0,\cdot)\right\|_{L^2}+\ln^{1/2}(3+t)\left\|\langle r\rangle^2 \Gamma^{\gamma+3}\partial_tU^{df}(0,\cdot)\right\|_{L^2},
    \end{aligned}
\end{equation}
which yields
\begin{equation}\label{668}
\begin{aligned}
H_2\lesssim&\frac{1}{\langle t\rangle^{3/2}}\chi_{[\mu/2]+4}(U^{cf})E_\mu^{1/2}(U^{cf})\left\|\Lambda U^{df}(0)\right\|_{\mu+1}\\
&+\frac{1}{\langle t\rangle^{3/2}}\chi_\mu(U^{cf})E_\mu^{1/2}(U^{cf})\left(E_{[\mu/2]+4}^{1/2}(U^{df})+\ln^{1/2}(3+t)\|\Lambda U^{df}(0)\|_{[\mu/2]+5}\right)\\
    &+\frac{1}{\langle t\rangle^{3/2}}\left(E_\mu^{1/2}(U^{cf})+\chi_\mu(U^{cf})\right)E_\mu^{1/2}(U^{df})M_\mu^{1/2}(U^{cf}).
\end{aligned}\end{equation}

For $N_1-N_4$, we benefit from the second null condition of the nonlinearities. Noting that
\begin{equation}
\langle t\rangle^{3/2}\lesssim \langle r\rangle\langle c_1t-r\rangle^{1/2}\langle c_2t-r\rangle^{1/2}\langle c_it-r\rangle^{1/2},\quad i=1\ \text{or}\ 2,
\end{equation}
we then have
\begin{equation}\label{6.60}
\begin{aligned}
N_1,N_2,N_3\lesssim &\frac{1}{\langle t\rangle^{3/2}}\left\|\langle r\rangle^{1/2}\langle c_1t-r\rangle^{1/2}\nabla\cdot \Gamma^{[\mu/3]+1} U^{cf}\right\|_{L^\infty}\left\|\langle r\rangle^{1/2}\partial \Gamma^{[\mu/3]+1} U\right\|_{L^\infty}\\
&\qquad\cdot\left\|\langle c_2t-r\rangle^{1/2}\nabla^\bot\cdot \Gamma^\mu U^{df}\right\|_{L^2}\left\|\langle c_1t-r\rangle^{1/2}\nabla\cdot \Gamma^\mu U^{cf}\right\|_{L^2}\\
&+\frac{1}{\langle t\rangle^{3/2}}\left\|\langle r\rangle^{1/2}\langle c_1t-r\rangle^{1/2}\nabla\cdot\Gamma^{[\mu/3]+1} U^{cf}\right\|_{L^\infty}\left\|\langle r\rangle^{1/2}\langle c_2t-r\rangle^{1/2}\nabla^\bot\cdot \Gamma^{[\mu/3]+1} U^{df}\right\|_{L^\infty}\\
&\qquad\cdot\left\|\nabla \Gamma^\mu U\right\|_{L^2}\left\|\langle c_1t-r\rangle^{1/2}\nabla\cdot \Gamma^\mu U^{cf}\right\|_{L^2}\\
&+\frac{1}{\langle t\rangle^{3/2}}\left\|\langle r\rangle^{1/2}\nabla\cdot\Gamma^{[\mu/3]+1} U^{cf}\right\|_{L^\infty}\left\|\langle r\rangle^{1/2}\langle c_2t-r\rangle^{1/2}\nabla^\bot\cdot \Gamma^{[\mu/3]+1} U^{df}\right\|_{L^\infty}\\
&\qquad\cdot\left\|\langle c_2t-r\rangle^{1/2}\nabla^\bot\cdot \Gamma^\mu U^{df}\right\|_{L^2}\left\|\langle c_1t-r\rangle^{1/2}\nabla\cdot \Gamma^\mu U^{cf}\right\|_{L^2}\\
&+\frac{1}{\langle t\rangle^{3/2}}\left\|\langle r\rangle^{1/2}\langle c_2t-r\rangle^{1/2}\nabla^\bot\cdot\Gamma^{[\mu/3]+1} U^{df}\right\|_{L^\infty}\left\|\langle r\rangle^{1/2}\langle c_2t-r\rangle^{1/2}\nabla^\bot\cdot \Gamma^{[\mu/3]+1} U^{df}\right\|_{L^\infty}\\
&\qquad\cdot\left\| \nabla\cdot \Gamma^\mu U^{cf}\right\|_{L^2}\left\|\langle c_1t-r\rangle^{1/2}\nabla\cdot \Gamma^\mu U^{cf}\right\|_{L^2}\\
\lesssim&\frac{1}{\langle t\rangle^{3/2}}\left(E_{[\mu/3]+4}(U)+E_{[\mu/3]+4}^{1/2}(U)\chi_{[\mu/3]+4}(U^{cf})\right)E_{\mu,\mu}^{1/2}(U^{cf},U^{df})M_{\mu}^{1/2}(U^{cf}).
\end{aligned}
\end{equation}
For the higher-order Morawetz type energy estimate, we need one more derivative for the divergence-free part $U^{df}$ in $N_4$, thus
\begin{equation}\label{6.61}
\begin{aligned}
N_4\lesssim& \frac{1}{\langle t\rangle^{3/2}}\left\|\langle r\rangle^{1/2}\langle c_2t-r\rangle^{1/2}\nabla^\bot\cdot\Gamma^{[\mu/3]+1} U^{df}\right\|_{L^\infty}\left\|\langle r\rangle^{1/2}\langle c_2t-r\rangle^{1/2}\nabla^\bot\cdot \Gamma^{[\mu/3]+1} U^{df}\right\|_{L^\infty}\\
&\qquad\cdot\left\| \nabla^\bot\cdot \Gamma^{\mu+1} U^{df}\right\|_{L^2}\left\|\langle c_1t-r\rangle^{1/2}\nabla\cdot \Gamma^\mu U^{cf}\right\|_{L^2}\\
\lesssim&\frac{1}{\langle t\rangle^{3/2}}E_{[\mu/3]+4}(U)E_{\mu+2}^{1/2}(U^{df})M_{\mu}^{1/2}(U^{cf}).
\end{aligned}
\end{equation}
Inserting \eqref{665}, \eqref{668}, \eqref{6.60} and \eqref{6.61} into \eqref{6.57}, we come to the lower-order energy estimate
\begin{equation}\label{6.62}
\begin{aligned}
  \frac{d \tilde{E}_{\mu}(U^{cf})}{dt}\leq&\frac{C}{\langle t\rangle^{3/2}}\chi_\kappa(U^{cf})E_\mu^{1/2}(U^{cf})\left(E_{[\mu/2]+4}^{1/2}(U^{df})+\ln^{1/2}(3+t)\|\Lambda U^{df}(0)\|_{\mu+1}\right)\\
    &+\frac{C}{\langle t\rangle^{3/2}}\left(E_\mu^{1/2}(U^{cf})+\chi_\mu(U^{cf})\right)E_{\mu+2}^{1/2}(U^{df})M_\mu^{1/2}(U^{cf}).
\end{aligned}\end{equation}

We are now in a position to give a detailed version of Theorem \ref{th61}.
\begin{thm}\label{thm62}
    Let $M_0>0$ and $0<\delta<\frac18$ be two given constants. For the Cauchy problem  \eqref{hweq} and \eqref{hwdata}, suppose that $U_0=U_0^{cf}+U_0^{df}$, $U_1=U_1^{cf}+U_1^{df}$ satisfy
    $$\nabla\cdot U_0^{df}=\nabla\cdot U_1^{df}=\nabla^\bot\cdot U_0^{cf}=\nabla^\bot\cdot U_1^{cf}=0$$
    and are supported in $\{x: |x|\leq 1^-\}$.
    Let $\kappa\ge14$ and $\mu=\kappa-2$. Suppose that 
    \begin{equation}
    \begin{aligned}
         M_{\kappa}^{1/2}\left(U^{cf}(0)\right)=&\sum_{|\alpha|\leq \kappa-1}\Bigg(\left\|\left(\frac{1+r}{8c_1}\right)^{1/2}\left(\Gamma_0^\alpha U_1^{cf}+c_1\partial_r\Gamma_0^\alpha U_0^{cf}\right)\right\|_{L^2}\\
         &\qquad+\left\|\left(\frac{1-r}{8c_1}\right)^{1/2}\left(\Gamma_0^\alpha U_1^{cf}-c_1\partial_r\Gamma_0^\alpha U_0^{cf}\right)\right\|_{L^2}\\
         &\qquad+\left\|\left(\frac{c_1}{2r^2}\right)^{1/2}\Omega\Gamma_0^\alpha U_0^{cf}\right\|_{L^2}+\left\|\left(\frac{c_1}{8}\right)^{1/2}\Gamma_0^\alpha U_0^{cf}\right\|_{L^2}\Bigg)\\
         \leq& M_0,
        \end{aligned}
        \end{equation}
        \begin{equation}\label{6.76}
        \begin{aligned}
          M_{\kappa+1}^{1/2}(U^{df}(0))=&\sum_{|\alpha|\leq \kappa-1}\Bigg(\sum_{a=0,1}\left\|\left(\frac{1+r}{8c_2}\right)^{1/2}\left(\nabla^a\Gamma_0^\alpha U_1^{df}+c_2\partial_r\nabla^a\Gamma_0^\alpha U_0^{df}\right)\right\|_{L^2}\\
         &\qquad+\sum_{a=0,1}\left\|\left(\frac{1-r}{8c_2}\right)^{1/2}\left(\nabla^a\Gamma_0^\alpha U_1^{df}-c_2\partial_r\nabla^a\Gamma_0^\alpha U_0^{df}\right)\right\|_{L^2}\\
         &\qquad+\sum_{a=0,1}\left\|\left(\frac{c_2}{2r^2}\right)^{1/2}\Omega\nabla^a\Gamma_0^\alpha U_0^{df}\right\|_{L^2}+\sum_{a=0,1}\left\|\left(\frac{c_2}{8}\right)^{1/2}\nabla^a\Gamma_0^\alpha U_0^{df}\right\|_{L^2}\Bigg)\\
         \leq &\varepsilon,
    \end{aligned} 
    \end{equation} 
\begin{equation}
    \begin{aligned}
        E_{\mu,\mu+1}^{1/2}(0,0)      =&\sum_{|\alpha|\leq\kappa-1}\left\|\Gamma_0^\alpha U_1^{cf}\right\|_{L^2}+c_1\left\|\nabla \cdot\Gamma_0^\alpha U_0^{cf}\right\|_{L^2}\\
        &+\sum_{|\alpha|\leq\kappa-1\atop a=0,1}\left\|\nabla^a\Gamma_0^\alpha U_1^{df}\right\|_{L^2}+c_2\left\|\nabla^a\nabla^\bot\cdot \Gamma_0^\alpha U_0^{df}\right\|_{L^2}\\
        \leq&\varepsilon,
    \end{aligned}
\end{equation}
\begin{equation}
    \begin{aligned}
        \|\Lambda U^{df}(0)\|_{\mu+1}=\sum_{\alpha=0}^{\mu}\left(\left\|\Gamma^\alpha U_0^{df}\right\|_{L^2}+\left\|\langle r\rangle^2 \Gamma^{\alpha}U^{df}_1\right\|_{L^2}\right)\leq \varepsilon.  
    \end{aligned}
\end{equation}
Then there exists a positive small constant $\epsilon_0$ which depends only on $\kappa$, $M_0$ and $\delta$, such that for any $0\leq \varepsilon\leq \epsilon_0$, 
\begin{itemize}
\item If $g_0=h_0=0$, the Cauchy problem  \eqref{hweq} and \eqref{hwdata} has an almost global classical solution satisfying
\begin{equation}
 M_{\kappa}^{1/2}(U^{cf})\leq C_0 M_0(1+t)^\delta;
\end{equation}
    \item If $g_0=g_2=h_0=0$, the Cauchy problem  \eqref{hweq} and \eqref{hwdata} has a unique global classical solution satisfying
\begin{equation}\label{678}
 M_{\kappa}^{1/2}(U^{cf})\leq C_0 M_0(1+t)^\delta,\qquad    E_{\mu}^{1/2}(U^{cf}) \leq \varepsilon A_0,
\end{equation}
\end{itemize}
where $C_0>1$ is a constant uniformly in $t$ and $A_0$ is a uniform constant  depending only on $C,C_0, M_0,\delta$.
\end{thm}
\emph{Proof.} We prove \eqref{678} by bootstrap argument. Recalling that $\tilde{E}_\mu(U^{cf}), \tilde{M}_\kappa(U^{cf})$ are equivalent to $E_\mu(U^{cf}), M_\kappa(U^{cf})$ respectively. For any $0<T<\infty$ and $0<t<T$, one may assume that
\begin{equation}\label{680}
\tilde{E}_{\mu}^{1/2}(U^{cf}) \leq \varepsilon A_0,\qquad \tilde{M}_{\kappa}^{1/2}(U^{cf})\leq C_0 M_0(1+t)^\delta.
\end{equation}
The identity \eqref{6.42} implies
$ M_{\kappa}(U^{df}(t))=M_{\kappa}(U^{df}(0))\leq\varepsilon,$ while
Lemma \ref{lem6.9} gives 
\begin{equation}
    \chi_\kappa(U^{cf})\lesssim E_\kappa^{1/2}(U^{cf})+\varepsilon M_{\kappa}^{1/2}(U^{cf})\lesssim M_{\kappa}^{1/2} ,
\end{equation}
\begin{equation}
    \chi_\mu(U^{cf})\lesssim E_\mu^{1/2}(U^{cf})+\varepsilon M_{\mu}^{1/2}(U^{cf}).
\end{equation}
From \eqref{6.62}, we then have
\begin{equation}\label{Emumu}
    \begin{aligned}
        \frac{d\tilde{E}_{\mu}(U^{cf})}{dt}\leq
        &\varepsilon^2 C A_0C_0M_0\langle t\rangle^{-3/2+\delta}\ln^{1/2}(3+t)+C C_0M_0\langle t\rangle^{-3/2+\delta}\tilde{E}_{\mu,\mu+1}(U^{cf},U^{df})\\
        &+\varepsilon^2 CA_0C_0^2M_0^2\langle t\rangle^{-3/2+2\delta}.
    \end{aligned}
\end{equation}
Let
\begin{equation}
    D=\int_0^TC C_0M_0\langle t\rangle^{-3/2+\delta}dt\leq \frac{2CC_0M_0}{1-2\delta}.
\end{equation}
Employing Gronwall's inequality to \eqref{el} yields
\begin{equation}\label{685}
    \begin{aligned}
     &\tilde{E}_{\mu}(U^{cf}(T))\\
     \leq&\varepsilon^2 e^D+ \varepsilon^2 C A_0C_0M_0e^D\int_0^Te^{-\int_0^tDd\tau}\left[\langle t\rangle^{-3/2+\delta}\ln^{1/2}(3+t)+C_0M_0\langle t\rangle^{-3/2+2\delta}\right]dt\\
     \leq&\varepsilon^2 e^{\frac{2CC_0M_0}{1-2\delta}}+\varepsilon^2 C A_0C_0M_0e^{\frac{4CC_0M_0}{1-2\delta}}B+\varepsilon^2 \frac{2}{1-4\delta}C A_0C_0^2M_0^2e^{\frac{4CC_0M_0}{1-2\delta}},
    \end{aligned}
\end{equation}
where $B=\int_0^T\langle t\rangle^{-3/2+\delta}\ln^{1/2}(3+t)dt$ can be bounded uniformly by some constant depending only on $\delta$. Choosing $$A_0>\max\left\{1,2e^{\frac{4CC_0M_0}{1-2\delta}}\left(1+CC_0M_0B+C_0^2M_0^2\frac{2}{1-4\delta}\right)\right\},$$
we get
\begin{equation}\label{687}
    \begin{aligned}
     \tilde{E}_{\mu}(U^{cf}(T))
     <\varepsilon^2A_0 e^{\frac{4CC_0M_0}{1-2\delta}}\left(1+ C C_0M_0B+ \frac{2}{1-4\delta}C C_0^2M_0^2\right)
     \leq \varepsilon^2A_0^2.
    \end{aligned}
\end{equation}

From \eqref{656}, we have
\begin{equation}
        \frac{d}{dt}\tilde{M}_{\kappa,\kappa+1}(U^{cf},U^{df})\leq \varepsilon CA_0\langle t\rangle^{-1}\tilde{M}_{\kappa,\kappa+1}(U^{cf}, U^{df}),
\end{equation}
which leads to
\begin{equation}\label{689}
        \tilde{M}_{\kappa,\kappa+1}(U^{cf},U^{df})(T)\leq M_0^2 (1+T)^{\epsilon_0 CA_0}< M_0^2 (1+T)^{2\delta},
\end{equation}
by choosing $\epsilon_0 (<\frac{2\delta}{CA_0})$ small enough.
Hence, we have improved the bootstrap assumption \eqref{680} via \eqref{687} and \eqref{689}. This completes the proof of Theorem \ref{thm62}.

\end{document}